\documentclass[aop,seceqn,MSNbibl,citesort,dvips]{arximspdf}
\usepackage{graphics}
\doi{10.1214/10-AOP545}
\volume{39}
\issue{1}
\pubyear{2011}
\firstpage{70}
\lastpage{103}

\makeatletter

\newcommand{\IP}{{\mathbb P}}
\newcommand{\IR}{{\mathbb R}}
\newcommand{\IN}{{\mathbb N}}
\newcommand{\IZ}{{\mathbb Z}}

\newcommand{\cA}{\mathcal{A}}
\newcommand{\cB}{\mathcal{B}}
\newcommand{\cT}{\mathcal{T}}
\newcommand{\cI}{\mathcal{I}}
\newcommand{\cH}{\mathcal{H}}
\newcommand{\cV}{\mathcal{V}}
\newcommand{\cW}{\mathcal{W}}

\newtheorem{theorem}{Theorem}[section]
\newtheorem{lemma}[theorem]{Lemma}
\newtheorem{lem}{Lemma}[section]
\newtheorem{corollary}[theorem]{Corollary}
\newtheorem{proposition}[theorem]{Proposition}
\newproclaim{remark}[theorem]{Remark}

\makeatother

\begin{document}
\begin{frontmatter}

\title{On the critical parameter of interlacement percolation in high
dimension}
\runtitle{Critical parameter of interlacement percolation}

\begin{aug}
\author[A]{\fnms{Alain-Sol} \snm{Sznitman}\corref{}\ead[label=e1]{sznitman@math.ethz.ch}}
\runauthor{A.-S. Sznitman}
\affiliation{ETH Z\"urich}
\address[A]{Departement Mathematik \\
ETH Z\"urich\\
CH-8092 Z\"urich\\
Switzerland\\
\printead{e1}} 
\end{aug}

\received{\smonth{1} \syear{2010}}
\revised{\smonth{3} \syear{2010}}

%
\begin{abstract}
The vacant set of random interlacements on $\IZ^d$, $d \ge3$, has
nontrivial percolative properties. It is known from Sznitman
[\textit{Ann. Math.} \textbf{171} (2010) 2039--2087], Sidoravicius and Sznitman
[\textit{Comm. Pure Appl. Math.} \textbf{62} (2009) 831--858] that
there is a nondegenerate critical value $u_*$ such that the vacant set
at level $u$ percolates when $u < u_*$ and does not percolate when $u >
u_*$. We derive here an asymptotic upper bound on $u_*$, as $d$ goes to
infinity, which complements the lower bound from Sznitman
[\textit{Probab. Theory Related Fields}, to appear]. Our main result shows that
$u_*$ is equivalent to $\log d$ for large $d$ and thus has the same
principal asymptotic behavior as the critical parameter attached to
random interlacements on $2d$-regular trees, which has been explicitly
computed in Teixeira [\textit{Electron. J. Probab.} \textbf{14}
(2009) 1604--1627].
\end{abstract}

%

\end{frontmatter}

\setcounter{section}{-1}
\section{Introduction}\label{sec0}

Random interlacements have proven useful in understanding how
trajectories of random walks can create large separating interfaces;
see \cite{Szni09c,Szni09d,CernTeixWind09}. In the
case of $\IZ^d$, $d \ge3$, it is known that the interlacement at
level $u \ge0$ is a random subset of $\IZ^d$, which is connected,
ergodic under translations
and infinite when $u$ is positive; see \cite{Szni07a}. The density of
this set monotonically increases from $0$ to $1$ as $u$ goes from $0$
to $\infty$. Its complement, the vacant set at level $u$, displays
nontrivial percolative properties. There is a critical value $u_*$ in
$(0,\infty)$ such that for $u < u_*$, the vacant set at level $u$ has
an infinite connected component which is unique (see
\cite{SidoSzni09a,Teix09a}) and, for $u > u_*$, only has finite
connected components; see \cite{Szni07a}. Little is known about $u_*$
and only recently was it shown that $u_*$ diverges when the dimension
$d$ tends to infinity; see \cite{Szni09e}. The aim of the present
article is to establish that $u_*$ is equivalent to $\log d$ as $d$
tends to infinity. In particular, this result shows that $u_*$ has the
same principal asymptotic behavior for large $d$ as the corresponding
critical parameter (which has been explicitly computed in
\cite{Teix09b}) attached to the percolation of the vacant set of random
interlacements on $2d$-regular trees.

We now describe the model. Precise definitions and pointers to the
literature appear in Section \ref{sec1}. Random interlacements are made of a
cloud of paths, which constitute a Poisson point process on the space
of doubly infinite $\IZ^d$-valued trajectories modulo time shift,
tending to infinity at positive and negative infinite times. The
nonnegative parameter $u$ mentioned above plays the role (roughly speaking)
of a multiplicative factor of the intensity measure of the Poisson
point process. Actually, one simultaneously constructs, on a suitable
probability space $(\Omega, \cA, \IP)$, the whole family $\cI^u$, $u
\ge0$, of random interlacements at level $u \ge0$ [cf. (\ref{1.30})].
They are the traces on $\IZ^d$ of the trajectories modulo time shift in
the cloud having labels at most $u$. The complement $\cV^u$ of $\cI^u$
in $\IZ^d$ is the vacant set at level $u$. It satisfies the following
identity:
%
%
\begin{equation}\label{0.1}
\IP[\cV^u \supseteq K] = \exp\{- u \operatorname{cap}(K)\} \qquad\mbox{for all
finite $K \subseteq\IZ^d$}.
\end{equation}
In fact, this formula provides a characterization of the law on $\{0,1\}
^{\IZ^d}$ of the indicator function of $\cV^u$ (cf. (2.16) of
\cite{Szni07a}). From Theorem 3.5 of \cite{Szni07a} and Theorem~3.4 of
\cite{SidoSzni09a}, one knows that there is a critical value $u_*$ in
$(0,\infty)$ such that:
%
%
\begin{eqnarray}\label{0.2}\hspace*{20pt}
&&\mbox{\hphantom{i}(i)} \quad \mbox{for $u > u_*$, $\IP$-a.s.}\qquad \mbox{all connected
components of $\cV^u$ are finite};
\nonumber\\
&&\mbox{(ii)} \quad \mbox{for $u < u_*$, $\IP$-a.s.}\qquad \mbox{there exists an
infinite connected}\\
&&\hspace*{123pt}\mbox{component in $\cV^u$}.\nonumber
\end{eqnarray}
From Theorem 0.1 of \cite{Szni09e}, one has the following asymptotic
lower bound on $u_*$ as $d$ tends to infinity:
%
%
\begin{equation}\label{0.3}
\liminf_d u_*\big/\log d \ge1 .
\end{equation}

The main aim of the present article is to show that the above lower
bound does capture the correct asymptotic behavior of $u_*$ and that
the following statement holds.
\begin{theorem}\label{theo0.1}
%
%
\begin{equation}\label{0.4}
\lim_d u_*\big/\log d= 1 .
\end{equation}
\end{theorem}

As a byproduct, this result shows that $u_*$ has the same principal
asymptotic behavior as the critical value attached to random
interlacements on $2d$-regular trees when $d$ goes to infinity; see
Proposition 5.2 of \cite{Teix09b}. We refer the reader to Remark~\ref
{rem4.1} for more on this matter. In addition, the proof of Theorem
\ref{theo0.1} also shows (cf. Remark \ref{rem4.1}) that
%
%
\begin{equation}\label{0.5}
\lim_d u_{**}\big/\log d = 1 ,
\end{equation}
where $u_{**} \in[u_*,\infty)$ is another critical value introduced
in \cite{Szni09c}. Informally, $u_{**}$ is the critical level above
which there is a polynomial decay in $L$ for the probability of
existence of a vacant crossing between a box of side length $L$ and the
complement of a concentric box of double side length. It is an
important and presently unresolved question whether $u_* = u_{**}$
actually holds. However, it is known that the connectivity function of
the vacant set at level $u$, that is, the probability that $0$ and a
(distant) $x$ are linked by a path in $\cV^u$ (i.e., the probability
of a vacant crossing at level $u$ between $0$ and $x$) has a stretched
exponential decay in $x$ when $u$ is bigger than $u_{**}$; see Theorem
0.1 of \cite{SidoSzni09b}.

We will briefly comment on the proof of Theorem \ref{theo0.1}. In view
of (\ref{0.3}), we only need to show that
%
%
\begin{equation}\label{0.6}
\limsup_d u_*\big/ \log d \le1 .
\end{equation}
As for Bernoulli bond or site percolation, similarities between what
happens on $\IZ^d$ and on $2d$-regular trees for large $d$ lurk in the
background of the proof. The statement corresponding to (\ref{0.6})
for Bernoulli percolation is an asymptotic lower bound for the critical
probability (a lower bound, not an upper bound, because the density of
$\cV^u$ decreases with $u$), whereas the required lower bound in the
Bernoulli percolation context follows from a short Peierls-type
argument (cf. \cite{BroaHamm57}, page 640, \cite{Kest90}, page 222,
or \cite{Grim99}, page 25); the proof of (\ref{0.6}) for random
interlacements is quite involved. The long-range dependence present in
the model is deeply felt.

An important feature of working in high dimension is that the $\ell
^1$-, Euclidean and $\ell^\infty$-distances all behave very
differently on $\IZ^d$; see (\ref{1.1}). At large enough scales
(i.e., Euclidean distance at least $d$), the Green function of the
simple random walk ``feels the invariance principle'' and is well
controlled by expressions of the type $(c \sqrt{d} / |\cdot|)^{d-2}$,
where $c$ does not depend on $d$ and \mbox{$|\cdot|$} stands for the
Euclidean norm; see Lemma \ref{lem1.1}. However, at shorter range, the
walk feels more of the tree-like nature of the space and the use of
bounds involving the $\ell^1$-distance becomes more pertinent
[cf. (\ref{1.14}) and Remark \ref{rem1.2}].

The above dichotomy permeates the proof of (\ref{0.6}). We use a
modification of the renormalization scheme (``for fixed $d$'') employed
in \cite{SidoSzni09b}. The renormalization scheme enables us to
transform certain local controls on the probability of vacant crossings
at level $u_0 = ( 1 + 5 \varepsilon) \log d$, $\varepsilon> 0$, small,
into controls on the probability of vacant crossings at arbitrary large
scales at a bigger level $u_\infty< ( 1 + 10 \varepsilon) \log d$.

The local estimates entering the initial step of the renormalization
scheme are developed in Section \ref{sec3}. These involve controls on the
existence of vacant crossings moving at $\ell^1$-distance
$c(\varepsilon)d$ from a box of side length $L_0 = d$ for the
interlacement at level $u_0$. The $2d$-regular tree model lurks behind
the control of these local crossings. The key estimates appear in
Theorem \ref{theo3.1} and Corollary \ref{cor3.4}. These estimates result from
an enhanced Peierls-type argument involving the consideration of what
happens in $\frac{c}{\varepsilon^2}$ $\ell^1$-balls, each having an
$\ell^1$-radius $c^\prime\varepsilon d$ and lying at mutual $\ell
^1$-distances of at least $c^{\prime\prime} d$. For this step, part
of the difficulty stems from the fact that the local estimates need to
be strong enough to overcome the combinatorial complexity involved in the
selection of the dyadic trees entering the renormalization scheme.

The renormalization scheme is developed in Section \ref{sec2}. It propagates
along an increasing sequence of levels $u_n$, with initial value $u_0 =
( 1 + 5 \varepsilon) \log d$ and limiting value $u_\infty< ( 1+ 10
\varepsilon) \log d$, uniform estimates on the probability of events
involving the presence of certain vacant crossings at level $u_n$.
Roughly speaking, these events correspond to the presence in $2^n$
boxes of side length $L_0(=d)$ of paths in~$\cV^{u_n}$. The boxes can
be thought of as the ``bottom leaves'' of a dyadic tree of depth $n$
and are well ``spread out'' within a box of side length $3L_n$, where
$L_n = \ell^n_0 L_0$ and $\ell_0 = d$. The paths start in each of
the $2^n$ boxes of side length $L_0$ and move at Euclidean (and hence
$\ell^1$-) distance of order $c(\varepsilon)d$ from the boxes. The
estimates are conducted uniformly over the possible dyadic trees
involved (cf. Propositions~\ref{prop2.1} and \ref{prop2.3}). The main
induction step in the above procedure (cf. Proposition \ref{prop2.1})
relies on the sprinkling technique introduced in \cite{Szni07a} to
control the long-range interactions. The rough idea is to introduce
more trajectories in the interlacement by letting the levels slightly
increase along the convergent sequence~$u_n$. In this way, one
dominates the long-range dependence induced by trajectories of the
interlacement traveling between distant boxes. In the present context,
the method uses, in an essential way, quantitative estimates on Harnack
constants in large Euclidean balls when the dimension $d$ goes to
infinity. These estimates crucially enter the proof of Proposition \ref
{prop2.3}. The bounds on the Harnack constants are derived in
Proposition \ref{prop1.3} with the help of the general Lem\-ma~\ref
{lemA.2} from the \hyperref[app]{Appendix}, which is an adaptation of Lemma 10.2 of
Grigoryan and Telcs~\cite{GrigTelc01}.

Let us now describe how this article is organized.

In Section \ref{sec1}, we introduce notation and recall several useful facts
concerning random walks and random interlacements. An important role is
played by the Green function bounds (see Lemma \ref{lem1.1}) and by
the bounds on Harnack constants; see Proposition \ref{prop1.3}.

In Section \ref{sec2}, we develop the renormalization scheme. It follows, with a
number of changes, the general line of \cite{SidoSzni09b}. The key
induction step appears in Proposition \ref{prop2.1}. The main
consequences of the renormalization scheme for the proof of Theorem
\ref{theo0.1} are stated in Proposition \ref{prop2.3}.

In Section \ref{sec3}, we derive the crucial local control on the existence of
vacant crossings at level $u_0$ traveling at $\ell^1$-distance of
order some suitable multiple of $d$. This local control is stated in
Theorem \ref{theo3.1}. It enables one to produce the required estimate
to initiate the renormalization scheme. This estimate can be found in
Corollary \ref{cor3.4}.

Section \ref{sec4} provides the proof of (\ref{0.6}). Combined with the lower
bound (\ref{0.3}) from~\cite{Szni09e}, this yields Theorem \ref
{theo0.1}. In Remark \ref{rem4.1}, we discuss some further questions
concerning the asymptotic behavior of $u_*$ for large $d$.

In the \hyperref[app]{Appendix}, we first derive, in Lemma \ref{lemA.1}, an elementary
inequality involved in proof of the Green function bounds from Lemma
\ref{lem1.1}. We then present, in Lemma \ref{lemA.2}, a general
result of independent interest providing controls on Harnack constants
in terms of killed Green functions
for general nearest-neighbor Markov chains on graphs.

Finally, let us explain the convention we use concerning constants.
Throughout the text, $c$ or $c^\prime$ denote positive constants with
values which can change from place to place. These constants
are independent of $d$. The numbered constants $c_0, c_1,\ldots$ are
fixed as the values of their first appearances in the text. Dependence
of constants on additional parameters appears in the notation, for
instance, $c(\varepsilon)$ denotes a constant depending on
$\varepsilon$.

\section{Notation and random walk estimates}\label{sec1}

In this section, we introduce further notation and gather various
useful estimates on simple random walk on $\IZ^d$ for large $d$.
Controls on the Green function and on Harnack constants in Euclidean
balls play an important role in the sequel. These can be found in Lemma
\ref{lem1.1} and Proposition \ref{prop1.3}. We also recall several
useful facts concerning random interlacements.

We let $\IN= \{0,1,2,\ldots\}$ denote the set of natural numbers.
Given a nonnegative real number $a$, we let $[a]$ denote the integer
part of $a$. We denote\vspace*{1pt} by $|\cdot|_1$, $|\cdot|$ and $|\cdot|_\infty
$ the $\ell^1$-, Euclidean and $\ell^\infty$-norms on $\IR^d$,
respectively. We have the following inequalities:
%
%
\begin{equation}\label{1.1}
|\cdot|_\infty\le| \cdot| \le| \cdot|_1 ,\qquad |\cdot| \le
\sqrt{d} |\cdot|_\infty,\qquad |\cdot|_1 \le\sqrt{d} |\cdot| .
\end{equation}
Unless explicitly stated otherwise, we tacitly assume that $d \ge3$.

By \textit{finite path}, we mean a sequence $x_0,\ldots, x_N$ in $\IZ
^d$, with $N \ge1$, which is such that $|x_{i+1} - x_i|_1 = 1$ for $0
\le i < N$. We sometimes write ``path'' in place of ``finite path''
when this causes no confusion. We denote by $B(x,r)$ and $S(x,r)$ the
closed ball and the closed sphere, respectively, with radius $r \ge0$
and center $x \in\IZ^d$. In the case of the $\ell^p$-distance where
$p=1$ or $\infty$, the corresponding objects are denoted by $B_p(x,r)$
and $S_p(x,r)$. For $A,B \subseteq\IZ^d$, we write $A + B$ for the
set of $x+y$ with $x$ in $A$ and $y$ in $B$, and $d(A,B) = \inf\{
|x-y|; x \in A, y \in B\}$ for the mutual Euclidean distance between
$A$ and $B$. We write $d_p(A,B)$, where $p=1$ or $\infty$, when the
$\ell^p$-distance is used instead. The notation $K \subset\subset\IZ
^d$ indicates that $K$ is a finite subset of $\IZ^d$. When $U$ is a
subset of $\IZ^d$, we write $|U|$ for the cardinality of $U$,
$\partial U = \{x \in U^c; \exists y \in U$, $|x-y|_1 = 1\}$ for the
boundary of $U$ and $\partial_{\mathrm{int}} U= \{x \in U$; $\exists y
\in U^c$, $|x-y|_1 = 1\}$ for the interior boundary of $U$. We also
write $\overline{U}$ in place of $U \cup\partial U$.

We denote by $W^+$ the set of nearest-neighbor $\IZ^d$-valued
trajectories defined for nonnegative times and tending to infinity. We
write $\cW_+$ and $X_n$, $n \ge0$, for the canonical $\sigma
$-algebra and the canonical process on $W_+$, respectively. We denote
by $\theta_n$, $n \ge0$, the canonical shift on $W_+$ so that $\theta
_n(w) = w(\cdot+ n)$ for $w \in W_+$ and $n \ge0$. Since $d \ge3$,
the simple random walk on $\IZ^d$ is transient and we write $P_x$ for
the restriction to the set $W_+$ of full measure of the canonical law
of the walk starting at $x \in\IZ^d$. When $\rho$ is a measure on
$\IZ^d$, we denote by $P_\rho$ the measure $\sum_{x \in\IZ^d} \rho
(x) P_x$ and by $E_\rho$ the corresponding expectation. Given $U
\subseteq\IZ^d$, we write $H_U = \inf\{n \ge0; X_n \in U\}$, $\widetilde
{H}_U = \inf\{n \ge1; X_n \in U\}$ and $T_U = \inf\{n \ge0; X_n
\notin U\}$ for the entrance time in $U$, the hitting time of $U$ and
the exit time from $U$, respectively. In the case of a singleton $\{x\}
$, we simply write $H_x$ and $\widetilde{H}_x$ for simplicity.

We let $g(\cdot,\cdot)$ stand for the Green function:
%
%
\begin{equation}\label{1.2}
g(x,x^\prime) = \sum_{n \ge0} P_x[X_n = x^\prime] \qquad\mbox{for
$x,x^\prime$ in $\IZ^d$}.
\end{equation}
The Green function is symmetric in its two variables and, due to
translation invariance, $g(x,x^\prime) = g(x^\prime- x) = g(x -
x^\prime)$, where
%
%
\begin{equation}\label{1.3}
g(x) = g(x,0) = g(0,x) \qquad\mbox{for } x \in\IZ^d .
\end{equation}
The $\ell^1$-distance is relevant for the description of the
short-range behavior of $g(\cdot)$ in high dimension [cf. Remark 1.3(1)
of \cite{Szni09e} and Remark \ref{rem1.2} below]; the Euclidean
distance becomes relevant in the description of the
``mid-to-long-range'' behavior of $g(\cdot)$. The following lemma will
be repeatedly used in the sequel. We recall that the convention
concerning constants is stated at the end of the \hyperref[sec0]{Introduction}.
\begin{lemma}\label{lem1.1}
%
%
\begin{eqnarray}
\label{1.4}
g(x) &\le& \bigl(c_0 \sqrt{d} / |x|\bigr)^{d-2} \qquad\mbox{for $|x| \ge d$}
\\
\label{1.5}
g(x) &\ge& \bigl(c_1 \sqrt{d} / |x|\bigr)^{d-2} \nonumber\\[-8pt]\\[-8pt]
\eqntext{\mbox{for $|x|^2 \ge d
|x|_\infty> 0$ (and, in particular, when $|x| \ge d$)}}
\\
\label{1.6}
\hspace*{22pt}P_x\bigl[H_{B(0,L)} < \infty\bigr] &\le& \biggl(\frac{c L}{|x|} \biggr)^{d-2}
\wedge1 \qquad\mbox{for $L \ge d, x \in\IZ^d$ (with $c \ge1$)}
\end{eqnarray}
\end{lemma}
\begin{pf}
We begin with the proof of (\ref{1.4}), (\ref{1.5}). To this end, we
denote by $p_t(u,v)$, $t \ge0$, $u,v \in\IZ$, the transition
probability of the simple random walk in continuous time on $\IZ$ with
exponential jumps of parameter $1$. The transition probability of the
simple random walk on $\IZ^d$ with exponential jumps of parameter $d$
can then be expressed as the product of one-dimensional transition
probabilities. Relating the continuous- and the discrete-time random
walks on $\IZ^d$, we thus find that
%
%
\begin{equation}\label{1.7}
g(x) = d \int^\infty_0 \prod^d_{i=1} p_t(0,x_i) \,dt\qquad
\mbox{for $x = (x_1,\ldots,x_d) \in\IZ^d$} .
\end{equation}
From Theorem 3.5 of \cite{Pang93} and the fact that the function
\[
F(\gamma) = - \log\bigl(\gamma+ \sqrt{\gamma^2 + 1} \bigr) +
\frac{1}{\gamma} \bigl(\sqrt{\gamma^2 + 1} - 1 \bigr),\qquad
\gamma> 0 ,
\]
appearing in Theorem 3.5 of \cite{Pang93} has derivative $-(1 + \sqrt
{\gamma^2 + 1})^{-1}$, tends to $0$ in $\gamma= 0$ and thus satisfies
the inequality $\log(1 + \frac{\gamma}{2}) \le- F(\gamma) \le\log
(1 + \gamma)$ for $\gamma\ge0$, we see that for suitable constants
$0 < \kappa< 1 < \kappa^\prime$, we have
%
%
\begin{eqnarray}\label{1.8}
&&\frac{1}{\kappa^\prime} (1 \vee t \vee
|u|)^{-{1/2}} \exp\biggl\{- |u| \log\biggl(1 + \kappa^\prime
\frac{|u|}{t} \biggr) \biggr\}\nonumber\\
&&\qquad \le p_t(0,u) \le
\frac{1}{\kappa} (1 \vee t \vee|u|)^{-{1/2}} \exp\biggl\{- |u| \log\biggl(1 + \kappa
\frac{|u|}{t} \biggr) \biggr\}\\
\eqntext{\mbox{for } t > 0, u \in\IZ.}
\end{eqnarray}
We now prove (\ref{1.4}) and thus assume that $|x| \ge d$. By (\ref
{1.7}), (\ref{1.8}), we bound $g(x)$ from above as follows (we also
use the inequality $d \le2^d$ and Lemma \ref{lemA.1} from the \hyperref[app]{Appendix}):
%
%
\begin{eqnarray}\label{1.9}
g(x) &\le& c^d \int^\infty_0 (1 \vee t)^{-{d/2}}
\exp\Biggl\{ - \sum^d_{i=1} |x_i| \log\biggl(1 + \kappa
\frac{|x_i|}{t} \biggr) \Biggr\} \,dt
\nonumber\\
&\stackrel{\mbox{\fontsize{8.36}{10.36}\selectfont{(\ref{A.1})}}}{\le} & c^d \int^\infty_0 (1 \vee
t)^{-{d/2}} \exp\biggl\{- |x| \log\biggl(1 + \kappa
\frac{|x|}{t} \biggr) \biggr\} \,dt
\nonumber\\[-8pt]\\[-8pt]
&\le& c^d \int_0^{\kappa|x|} (1 \vee t)^{-{d/2}}
\exp\biggl\{- |x| \log\biggl(1 + \kappa\frac
{|x|}{t} \biggr) \biggr\} \,dt
\nonumber\\
&&{} + c^d \int^\infty_{\kappa|x|} t^{-{d/2}} \exp
\biggl\{- \frac{\kappa|x|^2}{2 t} \biggr\} \,dt
,\nonumber
\end{eqnarray}
where, in the last step, we have used the inequality $\log(1 + \gamma
) \ge\frac{\gamma}{2}$ for $0 \le\gamma\le1$. Performing the
change of variable $s = \frac{\kappa|x|^2}{2t}$ in the last integral,
we see that the last term of (\ref{1.9}) is smaller than
%
%
\begin{equation}\label{1.10}\qquad
c^d |x|^{2-d} \int_0^{{|x|}/{2}} s^{{d}/{2} - 2}
e^{-s} \,d s \le c^d |x|^{2-d} \Gamma\biggl(\frac
{d}{2} - 1 \biggr) \le\bigl(c \sqrt{d} / |x|\bigr)^{d-2} ,
\end{equation}
using the asymptotic behavior of the gamma function in the last step
(cf. \cite{Olve74}, page~88).

As for the first integral in the last line of (\ref{1.9}), we note
that for $1 \le s \le\kappa|x|$, the function $s \rightarrow-\frac
{d}{2} \log s - |x| \log(1 + \kappa\frac{|x|}{s})$ has derivative
\[
-\frac{d}{2s} + \frac{|x|}{s}
\frac{\kappa|x|}{s + \kappa|x|} \stackrel{s \le
\kappa|x|}{\ge} - \frac{d}{2s} +
\frac{|x|}{2s} \stackrel{|x| \ge d}{\ge} 0
\]
and is hence nondecreasing. Thus, the first term in the last line of
(\ref{1.9}) is smaller than $c^d(\kappa|x|)^{-({d}/{2} - 1)} 2^{-|x|}$.

Observe that for $a \ge d$, $\frac{d-2}{2} \log a + a \log2 \ge
(d-2) \log\frac{a}{\sqrt{d}}$ (indeed, this inequality holds for
$a=d$ and $\frac{d-2}{2a}+ \log2 \ge\frac{d-2}{a}$ for $a \ge d$).
It follows that the first term in the last line of (\ref{1.9}) is at
most $(c \sqrt{d} / |x|)^{d-2}$. Together with (\ref{1.10}), this
completes the proof of (\ref{1.4}).

We now prove (\ref{1.5}) and assume that $x \not= 0$. Since $\log(1
+ \gamma) \le\gamma$ for $\gamma\ge0$, and $\kappa^\prime> 1$,
it follows from (\ref{1.7}), (\ref{1.8}) that
%
%
\begin{eqnarray}\label{1.11}
g(x) & \ge & c^d \int^\infty_{\kappa^\prime|x|_\infty}
t^{-{d/2}} \exp\biggl\{- \kappa^\prime
\frac{|x|}{t}^2 \biggr\} \,dt \nonumber\\
&\stackrel{s = {\kappa
^\prime
|x|^2}/{t}}{\ge}& c^d |x|^{2-d} \int_0^{
{|x|^2}/{|x|_\infty}}
s^{{d}/{2} -2} e^{-s} \,ds
\\
& \ge & \biggl(\frac{c\sqrt{d}}{|x|} \biggr)^{d-2}
\qquad\mbox{when $|x|^2 \ge d
|x|_\infty$} \nonumber
\end{eqnarray}
and (\ref{1.5}) follows.

Finally, (\ref{1.6}) is a routine consequence of the identity
%
%
\begin{equation}\label{1.12}\quad
g(x) = E_x \bigl[g\bigl(X_{H_{B(0,L)}}\bigr), H_{B(0,L)} < \infty\bigr] \qquad\mbox{for $L
\ge0$ and $x \in\IZ^d$} ,
\end{equation}
combined with (\ref{1.4}), (\ref{1.5}) and the fact that $\inf
_{B(0,L)} g \ge\inf_{\partial B(0,L)} g$.
\end{pf}
\begin{remark}\label{rem1.2}

(1) Although we will not need this fact in the sequel, let us mention
that the following lower bound complementing (\ref{1.6}) also holds:
%
%
\begin{equation}\label{1.13}\hspace*{25pt}
P_x\bigl[H_{B(0,L)} < \infty\bigr] \ge\biggl(\frac{c
L}{|x|} \biggr)^{d-2} \wedge1 \qquad\mbox{for $L \ge d$}, x \in
\IZ^d \mbox{ (with $c \le1$)} .
\end{equation}
Indeed, one uses (\ref{1.12}), together with (\ref{1.4}), (\ref
{1.5}), and, when $d + 1 \ge L( \ge d)$, the inequality $\sup
_{\partial_{\mathrm{int}} B(0,L)} g \le2d \sup_{\partial B (0,L)} g$,
which follows from the fact that $g$ is harmonic outside the origin
(the factor $2d$ can then be dominated by $\widetilde{c} {^{d-2}}$).

(2) Let us point out that when $x = ([d^\alpha],0,\ldots,0)$ with
$\frac{1}{2} < \alpha< 1$, the upper bound (\ref{1.4}) does not hold
when $d \ge c(\alpha)$. Indeed, it follows from (\ref{1.7}), (\ref
{1.8}) that
\[
g(x) \ge d \int^2_1 p_t(0,[d^\alpha]) p_t(0,0)^{d-1} \,dt
\stackrel{\mbox{\fontsize{8.36}{10.36}\selectfont{(\ref{1.8})}}}{\ge} c^d d^{-{\alpha/2}} \exp\{-
d^\alpha\log(1 + \kappa^\prime d^\alpha)\} ,
\]
which is much bigger than $(c_0 \sqrt{d} / |x|)^{d-2} \le c^{d-2} \exp
\{-(\alpha- \frac{1}{2}) (d-2) \log d\}$ for $d \ge c(\alpha)$.

(3) We recall from (\ref{1.11}) of \cite{Szni09e} that when $d \ge5$,
%
%
\begin{equation}\label{1.14}
g(x) \le\biggl(\frac{c_2 d}{|x|_1} \biggr)^{
{d/2} - 2} \qquad\mbox{for $x \in\IZ^d$} .
\end{equation}
The inequality is useful, for instance, when $|x| < d$, but $|x|_1 \ge
c_2 d$, a situation where (\ref{1.4}) is of no help. We will use (\ref
{1.14}) in Section \ref{sec3} when deriving local bounds on the connectivity
function of random interlacements at a level $u_0$ close to $\log d$;
see the proof of Theorem \ref{theo3.1}.

(4) The asymptotic behavior of $g(x)$ for $d$ fixed and large $x$ is
well known; see, for instance, \cite{HaraSlad92}, page 313, or
\cite{Lawl91}, page 31:
\[
\lim_{x \rightarrow\infty} \frac
{g(x)}{|x|^{d-2}} = \frac{d}{2} \Gamma
\biggl(\frac{d}{2} - 1 \biggr) \pi^{-{d/2}} .
\]
The asymptotic behavior of $g(\cdot)$ at the origin, or close to the
origin when $d$ tends to infinity, is also well known; see, for
instance, \cite{Mont56}, page 246, or \cite{Szni09e}, Remark~1.3(1).
On the other hand, the behavior of $g(\cdot)$ at intermediate scales
when $d$ tends to infinity seems much less well explored.
\end{remark}

The bounds on the Green function of Lemma \ref{lem1.1}, together with
Lemma \ref{lemA.2} from the \hyperref[app]{Appendix}, enable us to derive quantitative
controls on Harnack constants in suitably large Euclidean balls. These
bounds will be instrumental for the renormalization scheme developed in
the next section; see the proof of Lemma \ref{lem2.2}. First, we
recall some terminology. When $U \subseteq\IZ^d$, we say that a
function $u$ defined on $\overline{U}$ is harmonic in $U$
if, for all $x \in U$, $u(x) = \frac{1}{2d} \sum_{|e| = 1} u(x+e)$.
We can now state the following proposition.
\begin{proposition}[$(L \ge d)$]\label{prop1.3}
Setting $c_3 = 4 + 10 \frac{c_0}{c_1}$ [where $c_0 \ge c_1$---see
(\ref{1.4}), (\ref{1.5})], there exists $c > 1$ such that when $u$ is
a nonnegative function defined on $\overline{B(0,c_3 L)}$ and harmonic in
$B(0,c_3 L)$, we have
%
%
\begin{equation}\label{1.15}
\max_{B(0,L)} u \le c^d \min_{B(0,L)} u .
\end{equation}
\end{proposition}
\begin{pf}
We define $U_1 = B(0,L) \subseteq U_2 = B(0,4L) \subseteq U_3 = B(0,c_3
L)$. In view of Lemma \ref{lemA.2} from the \hyperref[app]{Appendix}, any $u$ as above
satisfies the inequality
\[
\max_{U_1} u \le K \min_{U_1} u ,
\]
where
%
%
\begin{equation}\label{1.16}
K = \max_{x,y \in U_1} \max_{z \in\partial_{\mathrm{int}} U_2} G_{U_3}(x,z) / G_{U_3}(y,z)
\end{equation}
and $G_{U_3}(\cdot,\cdot)$ stands for the Green function of the walk
killed outside $U_3$ [cf. (\ref{A.8})]. Applying the strong Markov
property at time $T_{U_3}$ and (\ref{1.2}), we obtain the following identity:
\[
G_{U_3}(y,z) = G(y,z) - E_y [G(X_{T_{U_3}}, z)] \qquad\mbox{for } y,z
\in\IZ^d .
\]
Hence, when $x,y \in U_1$ and $z \in\partial_{\mathrm{int}} U_2$, we see that
%
%
\begin{equation} \label{1.17}
G_{U_3}(x,z) \le G(x,z) \stackrel{\mbox{\fontsize{8.36}{10.36}\selectfont{(\ref{1.4})}}}{\le} \bigl(c_0 \sqrt{d}
/ (2L)\bigr)^{d-2}
\end{equation}
and
\begin{eqnarray}
\label{1.18}
G_{U_3}(y,z) & \ge & \bigl(c_1 \sqrt{d} / (5L)\bigr)^{d-2} - \bigl\{c_0 \sqrt{d} /
\bigl((c_3 - 4) L\bigr)\bigr\}^{d-2}
\nonumber\\
& = & \biggl(\frac{\sqrt{d}}{L} \biggr)^{d-2} \biggl(
\biggl(\frac{c_1}{5} \biggr)^{d-2} - \biggl(
\frac{c_0}{c_3 - 4} \biggr)^{d-2} \biggr)
\\
& = & \biggl(\frac{\sqrt{d}}{L} \biggr)^{d-2}
\biggl(\frac{c_1}{5} \biggr)^{d-2} \biggl(1 - \biggl(
\frac{1}{2}\biggr)^{d-2} \biggr) . \nonumber
\end{eqnarray}
We thus find that $K \le2 (\frac{5}{2} \frac{c_0}{c_1})^{d-2}$
and the claim (\ref{1.15}) follows.
\end{pf}

We now briefly review some notation and basic properties concerning the
equilibrium measure and the capacity. Given $K \subset\subset\IZ^d$,
we write $e_K$ for the equilibrium measure of $K$ and $\operatorname{cap}(K)$
for its total mass, the capacity of $K$:
%
%
\begin{eqnarray}\label{1.19}
e_K(x) &=& P_x[\widetilde{H}_K = \infty] 1_K(x),\qquad x \in\IZ^d,\nonumber\\[-8pt]\\[-8pt]
\operatorname{cap}(K) &=& \sum_{x \in K} P_x[\widetilde{H}_K = \infty]
.\nonumber
\end{eqnarray}
The capacity is subadditive [a straightforward consequence of (\ref{1.19})]:
%
%
\begin{equation}\label{1.20}
\operatorname{cap} (K \cup K^\prime) \le\operatorname{cap}(K) + \operatorname{cap}(K^\prime),
\qquad\mbox{for } K, K^\prime\subset\subset\IZ^d .
\end{equation}
One can also express the probability of entering $K$ in the following
well-known fashion:
%
%
\begin{equation}\label{1.21}
P_x[H_K < \infty] = \sum_{y \in K} g(x,y) e_K(y) \qquad\mbox{for
$x \in\IZ^d$} .
\end{equation}
Further, we have the following bound on the capacity of Euclidean balls:
%
%
\begin{equation}\label{1.22}
\operatorname{cap}(B(0,L)) \le\biggl(\frac{cL}{\sqrt{d}}
\biggr)^{d-2}\qquad
\mbox{for $L \ge d$} ,
\end{equation}
which follows from (\ref{1.5}), (\ref{1.6}) and (\ref{1.21}), by
letting $x$ tend to infinity.
\begin{remark}\label{rem1.4}
Although we will not need this
estimate in the sequel, let us mention that in a way analogous with
(\ref{1.4}), (\ref{1.13}) and (\ref{1.21}), one finds that
%
%
\begin{equation}\label{1.23}
\operatorname{cap} (B(0,L)) \ge\biggl(\frac{c L}{\sqrt{d}} \biggr)^{d-2}
\qquad\mbox{for $L \ge d$} .
\end{equation}
\end{remark}

We now turn to the description of random interlacements. We refer to
Section 1 of \cite{Szni07a} for details. We denote by $W$ the space of
doubly infinite nearest-neighbor $\IZ^d$-valued trajectories which
tend to infinity at positive and negative infinite times and by $W^*$
the space of equivalence classes of trajectories in $W$ modulo time
shift. We let $\pi^*$ stand for the canonical map from $W$ into $W^*$.
We write $\cW$ for the canonical $\sigma$-algebra on $W$ generated by
the canonical coordinates $X_n$, \mbox{$n \in\IZ$}, and $\cW^* = \{A
\subseteq W^*; (\pi^*)^{-1}(A) \in\cW\}$ for the largest $\sigma
$-algebra on $W^*$ for which $\pi^*\dvtx(W, \cW) \rightarrow(W^*, \cW
^*)$ is measurable. The canonical probability space for random
interlacements is now given as follows.

We consider the space of point measures on $W^* \times\IR_+$:
%
%
\begin{eqnarray}\label{1.24}\hspace*{22pt}
\Omega &=& \biggl\{\omega= \sum_{i \ge0} \delta_{(w^*_i,u_i)},
\mbox{with $(w^*_i,u_i) \in W^* \times\IR_+, i \ge0$ and } u \ge0,
\nonumber\\[-8pt]\\[-8pt]
&&\hspace*{47.27pt} w(W^*_K \times[0,u]) < \infty\mbox{ for any } K \subset\subset
\IZ^d\mbox{ and }u \ge0 \biggr\},\nonumber
\end{eqnarray}
where, for $K \subset\subset\IZ^d$, $W^*_K \subseteq W^*$ stands for
the set of trajectories modulo time shift that enter $K$, that is,
$W^*_K = \pi^*(W_K)$, where $W_K$ is the subset of $W$ of trajectories
that enter $K$.

We endow $\Omega$ with the $\sigma$-algebra $\cA$ generated by the
evaluation maps $\omega\rightarrow\omega(D)$, where $D$ runs over
the $\sigma$-algebra $\cW^* \times\cB(\IR_+)$, and with the
probability $\IP$ on $(\Omega, \cA)$, which is the Poisson measure
with intensity $\nu(d\omega^*) \,du$ giving finite mass to the sets
$W^*_K \times[0,u]$ for $K \subset\subset\IZ^d$, $u \ge0$, where
$\nu$ is the unique $\sigma$-finite measure on $(W^*,\cW^*)$ such
that for any $K \subset\subset\IZ^d$ (see Theorem 1.1 of \cite{Szni07a}),
%
%
\begin{equation}\label{1.25}
1_{W^*_K} \nu= \pi^* \circ Q_K ,
\end{equation}
$Q_K$ here denoting the finite measure on $W^0_K$, the subset of $W_K$
of trajectories which are for the first time in $K$ at time $0$ and
such that for $A,B \in\cW_+$ [we recall that $\cW_+$ is defined
above (\ref{1.2})] and $x \in\IZ^d$,
%
%
\begin{eqnarray}\label{1.26}
&&Q_K [(X_{-n})_{n \ge0} \in A, X_0 = x, (X_n)_{n \ge0} \in
B]\nonumber\\[-8pt]\\[-8pt]
&&\qquad =
P_x [A | \widetilde{H}_K = x] e_K(x) P_x[B] .\nonumber
\end{eqnarray}
For $K \subset\subset\IZ^d$, $u \ge0$, one defines on $(\Omega,
\cA)$ the following random variable valued in the set of finite point
measures on $(W_+, \cW_+)$:
%
%
\begin{eqnarray}\label{1.27}
\mu_{K,u}(dw) = \sum_{i \ge0} \delta_{(w_i^*)^{K,+}} 1\{w_i^* \in
W^*_K, u_i \le u\} \nonumber\\[-8pt]\\[-8pt]
\eqntext{\mbox{for } \displaystyle\omega= \sum_{i \ge0} \delta
_{(w_i^*,u_i)} \in\Omega,}
\end{eqnarray}
where, for $w^* \in W^*_K$, $(w^*)^{K,+}$ stands for the trajectory in
$W_+$ which follows $w^*$ step-by-step from the first time it enters $K$.

When $0 \le u^\prime< u$, one defines $\mu_{K,u^\prime,u}(dw)$ in an
analogous way to (\ref{1.27}), replacing the condition $u_i \le u$
with $u^\prime< u_i \le u$ in the right-hand side of (\ref{1.27}).
Then, for $0 \le u^\prime< u$, $K \subset\subset\IZ^d$, one finds that
%
%
\begin{equation}\label{1.28}
\begin{tabular}{p{260pt}}
$\mu_{K,u^\prime,u}$ and $\mu_{K,u^\prime}$ are independent
Poisson point processes
with respective intensity measures $(u-u^\prime) P_{e_K}$ and
$u^\prime P_{e_K}$.
\end{tabular}
\end{equation}
In addition, one has the identity
%
%
\begin{equation}\label{1.29}
\mu_{K,u} = \mu_{K,u^\prime} + \mu_{K,u^\prime, u} .
\end{equation}
Given $\omega\in\Omega$, the interlacement at level $u \ge0$ is the
following subset of $\IZ^d$:
%
%
\begin{eqnarray}\label{1.30}
\cI^u(\omega) & = & \bigcup_{u_i \le u} {\mathrm{range}}
(w_i^*) \hspace*{40pt}\nonumber\\[-8pt]\\[-8pt]
\eqntext{\mbox{if } \displaystyle\omega= \sum_{i \ge0} \delta_{(w_i^*,u_i)}
 =  \bigcup_{K \subset\subset\IZ^d} \bigcup_{w
\in\operatorname{Supp} \mu_{K,u}(\omega)}
w(\IN) ,}
\end{eqnarray}
where, for $w^* \in W^*$, range$(w^*) = w(\IN)$ for any $w \in W$,
with $\pi^*(w) = w^*$, and $\operatorname{Supp} \mu_{K,u}(\omega)$ refers
to the support of the point measure $\mu_{K,u}(\omega)$. The vacant
set at level $u$ is the complement of $\cI^u(\omega)$:
%
%
\begin{equation}\label{1.31}
\cV^u(\omega) = \IZ^d \setminus\cI^u(\omega)\qquad \mbox{for } u
\in\Omega, u \ge0 .
\end{equation}
One also has (cf. (1.54) of \cite{Szni07a})
%
%
\begin{equation}\label{1.32}\hspace*{28pt}
\cI^u(\omega) \cap K = \bigcup_{w \in\operatorname{Supp} \mu
_{K^\prime,u}(\omega)} w(\IN) \cap K \qquad\mbox{for } K \subset
K^\prime\subset\subset\IZ^d, u \ge0 .
\end{equation}
From (\ref{1.28}), one readily finds that, as mentioned in (\ref{0.1}),
%
%
\begin{equation}\label{1.33}
\IP[\cV^u \supseteq K] = \exp\{- u \operatorname{cap}(K)\} \qquad\mbox{for
all } K \subset\subset\IZ^d ,
\end{equation}
an identity that characterizes the law $Q_u$ on $\{0,1\}^{\IZ^d}$ of
the indicator function of $\cV^u(\omega)$; see also Remark 2.2(2) of
\cite{Szni07a}. This brings us to the conclusion of Section~\ref{sec1} and of
this brief review of some useful facts that we will use in the
following sections.

\section{From local to global: The renormalization scheme}\label{sec2}

In this section, we develop a renormalization scheme that follows, in
its broad lines, the strategy of~\cite{SidoSzni09b}. We introduce a
geometrically increasing sequence of length scales $L_n$, $n \ge0$,
and an increasing, but typically convergent, sequence of levels $u_n$,
$n \ge0$. When the sequence $u_n$ is sufficiently increasing
[cf. (\ref{2.19})], we are able to propagate from scale to scale
bounds on the key quantities $p_n(u_n)$ that appear in (\ref{2.17}).
Roughly speaking, these controls provide uniform upper bounds on the
probability that in a box at scale $L_n$, $2^n$ ``well spread'' boxes
at scale $L_0$ all witness certain crossing events at Euclidean
distance of order $c L_0$ in the vacant set at level~$u_n$.
Interactions are handled by the sprinkling technique originally
introduced in Section~3 of \cite{Szni07a}. The renormalization scheme
enables us to transform local estimates on the existence of vacant
crossings at scale $L_0$ in the vacant set at level $u_0$ into global
estimates on crossings at arbitrary scales in the vacant set at level
$u_\infty= \lim u_n$. The difficulty we encounter in the
implementation of the scheme stems from the fact that we want both
$u_0$ and $u_\infty$ to be ``slightly above'' the critical value
$u_*$; see (\ref{4.7}) and (\ref{0.3}). However, the local controls
on vacant crossings at level $u_0$, which we introduce into the
renormalization scheme and develop in the next section, require $L_0$
to be rather small, that is, of order $d$. We are then forced to keep a
tight control on the estimates we derive when $d$ goes to infinity. The
Green function and entrance probability estimates from Lemma \ref
{lem1.1}, together with the bounds on Harnack constants in Euclidean
balls from Proposition \ref{prop1.3}, play a pivotal role in this
scheme. The fact that the $\ell^\infty$- and Euclidean distances
behave very differently for large $d$ [see (\ref{1.1})] also forces
upon us some modifications of the geometric constructions in
\cite{SidoSzni09b}; see, for instance, (2.1) and (2.26). The main results of
this section are Proposition \ref{prop2.1}, which contains the main
induction step, and Proposition \ref{prop2.3}, which encapsulates the
estimates we will use in Section \ref{sec4}.

We consider the length scales
%
%
\begin{equation}\label{2.1}
L_0 \ge d,\qquad \widehat{L}_0 = \bigl(\sqrt{d} + R\bigr) L_0\qquad \mbox{with $R \ge1$}
\end{equation}
as well as
%
%
\begin{eqnarray}\label{2.2}
L_n = \ell_0^n L_0\hspace*{140pt}\nonumber\\[2pt]\\[-22pt]
\eqntext{\begin{tabular}{p{238pt}}
for $n \ge1$, where $\ell_0 \ge
1000 \dfrac{c_0}{c_1} \bigl(\sqrt{d} + R\bigr)$ is an integer
multiple of $100$ (we recall that $c_0 \ge c_1$; cf. Lemma
\ref{lem1.1}).\hspace*{-18pt}
\end{tabular}}
\end{eqnarray}
We organize $\IZ^d$ in a hierarchical way, with $L_0$ being the finest
scale and $L_1 < L_2 < \cdots$ being coarser and coarser scales.
Crossing events at the finest scale will involve the length scale
$\widehat
{L}_0$. We introduce the following set of labels of boxes at level $n
\ge0$:
%
%
\begin{equation}\label{2.3}
I_n = \{n\} \times\IZ^d .
\end{equation}
To each $m = (n,i) \in I_n$, $n \ge0$, we attach the box
%
%
\begin{equation}\label{2.4}
C_m = \bigl(i L_n + [0,L_n)^d\bigr) \cap\IZ^d .
\end{equation}
In addition, when $n \ge1$, we define
%
%
\begin{equation}\label{2.5}
\widetilde{C}_m = \bigcup_{m^\prime\in I_n, d_\infty(C_{m^\prime
}, C_m) \le1} C_{m^\prime} ( \mbox{$\supseteq$} C_m) .
\end{equation}
On the other hand, when $n = 0$ and $m = (0,i) \in I_0$, we define instead
%
%
\begin{equation}\label{2.6}
\widetilde{C}_m = B(i L_0, \widehat{L}_0) \stackrel{\mbox{\fontsize{8.36}{10.36}\selectfont{(\ref{1.1}),
(\ref{2.1})}}}{\supseteq} \bigcup_{x \in C_m} B(x,R L_0) ( \mbox{$\supseteq$}
C_m) .
\end{equation}
The above definitions slightly differ from (2.3) in \cite{SidoSzni09b}
due to the special treatment of the scale $n=0$. It is relevant here to
use Euclidean balls and, thanks to (\ref{1.6}) of Lemma \ref{lem1.1},
to have a good control on the entrance probability of a simple random
walk in $\widetilde{C}_m$. The radius of these balls has to be chosen
sufficiently large so that we can show that crossing events at the
bottom scale, from $C_m$ to $\partial_{\mathrm{int}} \widetilde{C}_m$, are
unlikely (this will be done in the next section).

\begin{figure}

\includegraphics{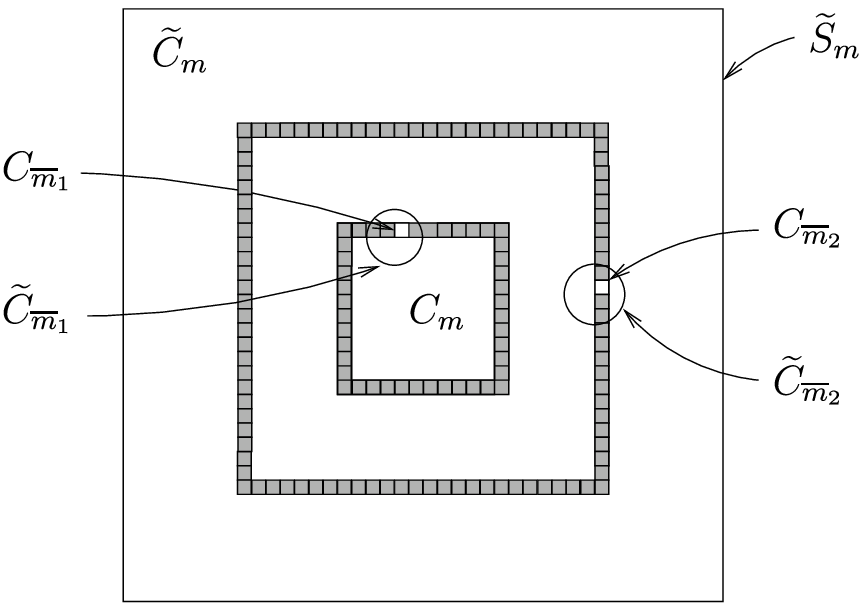}

\caption{An illustration of the boxes $C_{\overline{m}_i}$ and balls
$\widetilde
{C}_{\overline{m}_i}$, $i=1,2$, when $m$ belongs to $I_1$.}
\label{fig1}
\end{figure}

We then write $S_m = \partial_{\mathrm{int}} C_m $ and $\widetilde{S}_m =
\partial_{\mathrm{int}} \widetilde{C}_m$ for $m \in I_n$, $n \ge0$. Given $m
\in
I_n$ with $n \ge1$, we consider $\cH_1(m)$, $\cH_2(m) \subseteq
I_{n-1}$ defined by
%
%
\begin{eqnarray}\label{2.7}
\cH_1(m) & = & \{\overline{m} \in I_{n-1}; C_{\overline{m}} \subseteq C_m
\mbox{ and } C_{\overline{m}} \cap S_m \not= \varnothing\},
\nonumber\\[-8pt]\\[-8pt]
\cH_2(m) & = & \biggl\{\overline{m} \in I_{n-1}; C_{\overline{m}} \cap\biggl\{ z
\in\IZ^d; d_\infty(z,C_m) = \frac{L_n}{2}
\biggr\}\not= \varnothing
\biggr\} .\nonumber
\end{eqnarray}
We thus see that for $n \ge1$, $m \in I_n$, one has (see also
Figure \ref{fig1}):
%
%
\begin{eqnarray}\label{2.8}\qquad
&&\overline{m}_1 \in \cH_1(m),\qquad
\overline{m}_2 \in \cH_2(m)\nonumber\\[-8pt]\\[-8pt]
&&\mbox{implies that }
\widetilde{C}_{\overline{m}_1} \cap\widetilde{C}_{\overline{m}_2} =
\varnothing\mbox{ and }
\widetilde{C}_{\overline{m}_1} \cup\widetilde{C}_{\overline{m}_2}
\subseteq\widetilde{C}_m
\nonumber
\end{eqnarray}
[in the case $n=1$, we use the lower bound on $\ell_0$ in (\ref{2.2})
as well as (\ref{1.1})].

%
%

Then, to each $m \in I_n$, $n \ge0$, we associate a collection
$\Lambda_m$ of ``binary trees of depth $n$.'' More precisely, we
define $\Lambda_m$ to be the collection of subsets $\cT$ of $\bigcup
_{0 \le k \le n} I_k$ such that, writing $\cT^k = \cT\cap I_k$, we have
%
%
\begin{eqnarray}\label{2.9}\hspace*{29pt}
&\cT^n = \{m\},&
\\
%
%
\label{2.10}
&
\begin{tabular}{p{310pt}}
any $m^\prime\in\cT^k$, $1 \le k \le n$, has two
``descendants,'' $\overline{m}_i(m^\prime) \in\cH_i(m^\prime)$, \mbox{$i = 1,2$},
such that $\cT^{k-1} = \bigcup_{m^\prime\in\cT^k} \{
\overline{m}_1 (m^\prime), \overline{m}_2(m^\prime)\}$.
\end{tabular}
&
\end{eqnarray}
For each $\cT\in\Lambda_m$ and $m^\prime\in\cT$, one can then
define the subtree of ``descendants of $m^\prime$ in $\cT$'' via
%
%
\begin{equation}\label{2.11}
\cT_{m^\prime} = \{m^{\prime\prime} \in\cT; \widetilde{C}_{m^{\prime
\prime}} \subseteq\widetilde{C}_{m^\prime}\} (\mbox{$\in$}\Lambda_{m^\prime}) .
\end{equation}
Given $1 \le k \le n$, $m^\prime\in\cT^k$, one thus has the
following partition of $\cT_{m^\prime}$:
%
%
\begin{equation}\label{2.12}
\cT_{m^\prime} = \{m^\prime\} \cup\cT_{\overline{m}_1(m^\prime)} \cup
\cT_{\overline{m}_2(m^\prime)} .
\end{equation}
In addition, we have the following rough bound on the collection
$\Lambda_m$ of binary trees attached to $m \in I_n$:
%
%
\begin{equation}\label{2.13}\hspace*{28pt}
|\Lambda_m| \le(c_4 \ell_0)^{2(d-1)} (c_4 \ell_0)^{4(d-1)}
\cdots(c_4 \ell_0)^{2^n(d-1)} = (c_4 \ell_0)^{2(d-1) (2^n-1)} ,
\end{equation}
where we have used the rough bound for $m^\prime\in I_k$, $1 \le k \le
n$, and, for $i = 1,2$,
\[
|\cH_i(m^\prime)| \le2d \biggl(c \frac
{L_k}{L_{k-1}} \biggr)^{d-1} = 2d (c \ell_0)^{d-1} \le(c_4 \ell
_0)^{d-1} \qquad\mbox{for some $c_4 > 1$} .
\]
We then introduce, for $u \ge0$, $m \in I_n$, with $n \ge0$, the event
%
%
\begin{equation}\label{2.14}
A^u_m = \bigl\{C_m \stackrel{\cV^u}{\longleftrightarrow} \widetilde
{S}_m \bigr\} ,
\end{equation}
where the expression in the right-hand side of (\ref{2.14}) denotes
the collection of $\omega$ in $\Omega$ such that there is a path
between $C_m$ and $\widetilde{S}_m$ contained in $\cV^u$. In an analogous
fashion to Lemma 2.1 of \cite{SidoSzni09b}, $A^u_m$ ``cascades down to
the bottom scale'' because any path originating in $C_m$ and ending in
$\widetilde{S}_m$ must go through some $C_{\overline{m}_1}$, $\overline
{m}_1 \in\cH
_1(m)$, reach $\widetilde{S}_{\overline{m}_1}$ and then go through some
$C_{\overline
{m}_2}, \overline{m}_2\in\cH_2(m)$, and reach $\widetilde{S}_{\overline
{m}_2}$. Thus,
similarly to Lemma 2.1 of \cite{SidoSzni09b}, we find that defining
for $u \ge0$, $n \ge0$, $m \in I_n$ and $\cT\in\Lambda_m$
%
%
\begin{equation}\label{2.15}
A^u_\cT= \bigcap_{m^\prime\in\cT^0} A^u_{m^\prime}
\qquad\mbox{(recall that $\cT^0 = \cT\cap I_0$)} ,
\end{equation}
one has the inclusion
%
%
\begin{equation}\label{2.16}
A^u_m \subseteq\bigcup_{\cT\in\Lambda_m} A^u_{\cT} .
\end{equation}
We then introduce the key quantity
%
%
\begin{equation}\label{2.17}
p_n(u) = \sup_{\cT\in\Lambda_m} \IP[A^u_\cT],\qquad u \ge0,
n \ge0 \mbox{ with $m \in I_n$ arbitrary} ,
\end{equation}
which is well defined due to translation invariance, and find that
%
%
\begin{equation}\label{2.18}
\IP[A^u_m] \le|\Lambda_m| p_n(u) \qquad\mbox{for $u \ge0, n \ge0$} .
\end{equation}
The heart of the matter is now to find a recurrence relation bounding
$p_{n+1}(u_{n+1})$ in terms of $p_n(u_n)$ for suitably increasing
sequences $u_n$ (we are actually interested in increasing, but
convergent, sequences). We recall that $R$ appears in (\ref{2.1}).
\begin{proposition}\label{prop2.1}
There exist positive constants $c_5,c_6,c$ such that if $\ell_0 \ge
c(\sqrt{d} + R)$, then, for any increasing sequences $u_n$, $n \ge0$,
in $(0,\infty)$ and nondecreasing sequences $r_n$, $n \ge0$, of
positive integers such that
%
%
\begin{equation}\label{2.19}
u_{n+1} \ge u_n \biggl(1 + \frac{\widehat{L}_0}{L_0}
\biggl(\frac{c_5}{\ell_0} \biggr)^{(n+1) (d-2)}
\biggr)^{r_n+1} \qquad\mbox{for all $n \ge0$} ,
\end{equation}
one has, for all $n \ge0$,
%
%
\begin{eqnarray}\label{2.20}
p_{n+1}(u_{n+1}) &\le& p_n(u_{n+1})\biggl(p_n(u_n) + u_n \biggl(
\frac{\widehat{L}_0}{\sqrt{d}} \biggr)^{(d-2)}\nonumber\\[-8pt]\\[-8pt]
&&\hspace*{45.17pt}{}\times  \biggl(4^n \biggl(c_6
\frac{\widehat{L}_0}{L_0} \biggr)^{(d-2)} \ell
_0^{-(n+1)(d-2)} \biggr)^{r_n} \biggr)\nonumber
\end{eqnarray}
[note that $p_n(\cdot)$ is nonincreasing so that $p_n(u_{n+1}) \le p_n(u_n)$].
\end{proposition}
\begin{pf}
The proof of Proposition \ref{prop2.1} is an adaptation of the proof
of Proposition 2.2 of \cite{SidoSzni09b}, which will be sketched below
with some modifications which we will highlight.

One considers some $n \ge0$, $m \in I_{n+1}$, $\cT\in\Lambda_m$ and
writes $\overline{m}_1, \overline{m}_2$ for the unique elements of $\cH_1(m)$,
$\cH_2(m)$ in $\cT^n$ ($=\cT\cap I_n$). One also
writes $u^\prime$ and~$u$, with $0 < u^\prime< u$, in place of $u_n$
and $u_{n+1}$.

If $\overline{\cT} \in\Lambda_{\overline{m}}$ with $\overline{m} \in
I_n$, then one
defines, for $\mu$, a point process on $W_+$ defined on $\Omega$
(i.e., a measurable map from $\Omega$ into the space of point measures
on $W_+$):
%
%
\begin{eqnarray}\label{2.21}\qquad
A_{\overline{\cT}}(\mu) &=& \bigcap_{m^\prime\in\overline{\cT} \cap
I_0} \biggl\{  \omega\in\Omega;
\mbox{there is a path in }
\nonumber\\[-8pt]\\[-8pt]
&&{}\hspace*{39.7pt}\widetilde{C}_{m^\prime} \Bigm\backslash\biggl(\bigcup
_{w \in\operatorname{Supp} \mu(\omega)} w(\IN) \biggr)\mbox{ joining }C_{m^\prime}\mbox{ with }
\widetilde{S}_{m^\prime} \biggr\}.\nonumber
\end{eqnarray}
As in (\ref{2.19}) of \cite{SidoSzni09b}, using independence, we have
the bound
%
%
\begin{equation}\label{2.22}
\IP[A^u_{\cT}] \le p_n(u) \IP[A_{\overline{\cT}_2} (\mu_{2,2})] ,
\end{equation}
where $\overline{\cT}_2$ stands for $\cT_{\overline{m}_2}$ and we have
decomposed the point process $\mu_{V,u}$ [see (\ref{1.27})], where
%
%
\begin{eqnarray} \label{2.23}
\hspace*{37.4pt}\hspace*{67.4pt} V & = & \widehat{C}_1 \cup\widehat{C}_2
\\
\label{2.24}
&&\mbox{with }\widehat{C}_i = \bigcup_{m^\prime\in\overline{\cT}_i \cap I_0}
\widetilde{C}_{m^\prime} \subseteq\widetilde{C}_{\overline{m}_i}
\qquad\mbox{for $i=1,2$}
\end{eqnarray}
(i.e., a union of $2^n$ pairwise disjoint Euclidean balls of radius
$\widehat{L}_0$), into a sum of independent Poisson processes via
%
%
\begin{equation}\label{2.25}
\mu_{V,u} = \mu_{1,1} + \mu_{1,2} + \mu_{2,1} + \mu_{2,2} ,
\end{equation}
where, for $i \not= j$ in $\{1,2\}$, we have set
\[
\mu_{i,j} = 1 \{X_0 \in\widehat{C}_i, H_{\widehat{C}_j} < \infty\} \mu
_{V,u} \quad\mbox{and}\quad \mu_{i,i} = 1\{X_0 \in\widehat{C}_i, H_{\widehat
{C}_j} = \infty\} \mu_{V,u} .
\]
One introduces similar decompositions for $\mu_{V,u^\prime}$ in terms
of analogously defined point processes $\mu^\prime_{i,j}$, $1 \le i,j
\le2$, and for $\mu_{V,u^\prime,u}$ in terms of $\mu^*_{i,j}$, $1
\le i,j \le2$.

The heart of the matter is to bound $\IP[A_{\overline{\cT}_2}(\mu_{2,2})]
= \IP[A_{\overline{\cT}_2}(\mu^\prime_{2,2} + \mu^*_{2,2})]$, which
appears in the right-hand side of (\ref{2.22}), in terms of
$p_n(u^\prime)$ when $u-u^\prime$ is not too small. For this purpose,
we employ the sprinkling technique of \cite{Szni07a} and, loosely
speaking, establish that $\mu^*_{2,2}$ dominates ``up to small
corrections'' the contribution of $\mu^\prime_{2,1} + \mu^\prime
_{1,2}$ in $\IP[A^{u^\prime}_{\overline{\cT}_2}] = \IP[A_{\overline{\cT
}_2}(\mu^\prime_{2,2} + \mu^\prime_{2,1} + \mu^\prime_{1,2})]$.

With this in mind, we define a neighborhood $U$ of $\widetilde
{C}_{\overline
{m}_2}$ (and, in contrast to (2.20) of \cite{SidoSzni09b}, we do not
define $U$ as the $\ell^\infty$-neighborhood of $\widetilde
{C}_{\overline{m}_2}$
of size $\frac{L_{n+1}}{10}$). Instead, if $\overline{m}_2 =
(n,\overline{i}_2)
\in I_n$ [see (\ref{2.3})], we define $U$ as the following Euclidean
ball (which is much smaller than the corresponding $\ell^\infty$-ball
of same radius):
%
%
\begin{equation}\label{2.26}
\quad U = B \biggl(\overline{i}_2 L_n, \frac{L_{n+1}}{10} \biggr)
\supseteq\widetilde{C}_{\overline{m}_2}\qquad \mbox{using (\ref{2.1}), (\ref
{2.2}), (\ref{2.5}), (\ref{2.6})} .
\end{equation}
We then have the following important controls on Euclidean distances:
%
%
\begin{eqnarray}\label{2.27}
d(\partial U, \widehat{C}_2) & \ge & \frac{L_{n+1}}{10} -
3 \sqrt{d} L_n \stackrel{\mbox{\fontsize{8.36}{10.36}\selectfont{(\ref{2.1}), (\ref{2.2})}}}{>}
\frac{L_{n+1}}{20} \qquad\mbox{when $n \ge1$}
\nonumber\\[-8pt]\\[-8pt]
& \ge & \frac{L_{n+1}}{10} - \widehat{L}_0 \stackrel
{\mbox{\fontsize{8.36}{10.36}\selectfont{(\ref{2.1}), (\ref{2.2})}}}{>} \frac{L_{n+1}}{20}
\qquad\mbox{when $n =
0$}, \nonumber
\end{eqnarray}
and we have used in the first line the fact that $\widehat{C}_2
\subseteq
\widetilde{C}_{\overline{m}_2}$ when $m \ge1$; see (\ref{2.8}). Using similar
considerations, we find that
%
%
\begin{eqnarray}\label{2.28}\hspace*{28pt}
d(\partial U, \widehat{C}_1) & \ge & \frac{L_{n+1}}{2} -
L_n - L_n - \frac{L_{n+1}}{10} - 1 >
\frac{L_{n+1}}{20} \qquad\mbox{when $n \ge1$}
\nonumber\\[-8pt]\\[-8pt]
& \ge & \frac{L_{n+1}}{2} - \widehat{L}_0 - L_0 -
\frac{L_{n+1}}{10} - 1 > \frac
{L_{n+1}}{20} \qquad\mbox{when $n=
0$} .\nonumber
\end{eqnarray}
Since $V = \widehat{C}_1 \cup\widehat{C}_2$, we have established that
%
%
\begin{equation}\label{2.29}
d(\partial U, V) > \frac{L_{n+1}}{20} .
\end{equation}
We then introduce the successive times of return to $\widehat{C}_2$ and
departure from $U$:
%
%
\begin{eqnarray}\label{2.30}\hspace*{38pt}
R_1 & = & H_{\widehat{C}_2},\qquad D_1 = T_U \circ\theta_{R_1} + R_1 \quad\mbox
{and}\quad\mbox{for $k \ge1$, by induction,}
\\
R_{k+1} & = & R_1 \circ\theta_{D_k} + D_k,\qquad D_{k+1} = D_1 \circ\theta
_{D_k} +
D_k, \nonumber
\end{eqnarray}
so that $0 \le R_1 \le D_1 \le\cdots\le R_k \le D_k \le\cdots\le
\infty$.

Letting $r \ge1$ play the role of $r_n$ in (\ref{2.19}), (\ref
{2.20}), we further introduce the decompositions
%
%
\begin{eqnarray}\label{2.31}
\mu^\prime_{2,1} & = & \sum_{1 \le\ell\le r} \rho^\ell_{2,1} +
\overline{\rho}_{2,1},\qquad \mu^\prime_{1,2} = \sum_{1 \le\ell\le r}
\rho^\ell_{1,2} + \overline{\rho}_{1,2} ,
\nonumber\\[-8pt]\\[-8pt]
\mu^*_{2,2} & = & \sum_{1 \le\ell\le r} \rho^\ell_{2,2} + \overline
{\rho}_{2,2}
,\nonumber
\end{eqnarray}
where, for $i \not= j$ in $\{1,2\}$ and $\ell\ge1$, we have set
\begin{eqnarray*}
\rho^\ell_{i,j} & = & 1\{R_\ell< D_\ell< R_{\ell+ 1} = \infty\}
\mu^\prime_{i,j},\\
\overline{\rho}_{i,j} & = & 1\{R_{r + 1} < \infty\} \mu
^\prime_{i,j},
\\
\rho^\ell_{2,2} & = & 1\{R_\ell< D_\ell< R_{\ell+ 1} = \infty\}
\mu^*_{2,2}
\end{eqnarray*}
and
\[
\overline{\rho}_{2,2} = 1\{R_{r + 1} < \infty
\} \mu^*_{2,2} .
\]
The point processes $\overline{\rho}_{1,2}$ and $\overline{\rho}_{2,2}$ play
the role of correction terms, eventually responsible for the last term
in the right-hand side of (\ref{2.20}). The bounds we derive on the
intensity measures $\overline{\xi}_{2,1}$ and $\overline{\xi}_{1,2}$ of
$\overline
{\rho}_{2,1}$ and $\overline{\rho}_{1,2}$ depart from
(2.26), (2.27) in \cite{SidoSzni09b}. We write
%
%
\begin{eqnarray}\label{2.32}
\overline{\xi}_{2,1}(W_+) &=& u^\prime
P_{e_V} [X_0 \in\widehat{C}_2, H_{\widehat{C}_1} < \infty, R_{r+1} <
\infty]
\nonumber\\
& \stackrel{\mbox{\fontsize{8.36}{10.36}\selectfont{(\ref{1.19})}}}{\le} &
u^\prime\operatorname{cap}(\widehat{C}_2) \sup_{x \in\widehat{C}_2} P_x
[R_{r+1} < \infty]
\\
& \stackrel{\mathrm{strong}\ \mathrm{Markov}}{\le} &
u^\prime\operatorname{cap}(\widehat{C}_2) \Bigl( \sup_{x \in\partial U}
P_x[H_{\widehat{C}_2} < \infty]\Bigr)^r .\nonumber
\end{eqnarray}
Combining (\ref{1.6}) and (\ref{2.29}), we find that
%
%
\begin{equation}\label{2.33}
\sup_{x \in\partial U} P_x [H_{\widehat{C}_2} < \infty] \le
2^n \biggl( c \frac{\widehat{L}_0}{L_{n+1}}
\biggr)^{(d-2)} \stackrel{\mbox{\fontsize{8.36}{10.36}\selectfont{(\ref{2.1}), (\ref{2.2})}}}{=} 2^n \biggl(c
\frac{\widehat{L}_0}{L_0} \ell_0^{-(n+1)} \biggr)^{(d-2)} .\hspace*{-35pt}
\end{equation}
Moreover, from (\ref{1.20}), (\ref{1.22}), we have
%
%
\begin{equation}\label{2.34}
\operatorname{cap}(\widehat{C}_2) \le2^n \biggl(c \frac{\widehat
{L}_0}{\sqrt{d}} \biggr)^{(d-2)}
\end{equation}
and hence
%
%
\begin{equation}\label{2.35}
\overline{\xi_{2,1}}(W_+) \le u^\prime\biggl(\frac{\widehat
{L}_0}{\sqrt{d}} \biggr)^{(d-2)} \biggl(4^n \biggl(c
\frac{\widehat{L}_0}{L_0} \biggr)^{(d-2)} \ell_0^{-(n+1)(d-2)} \biggr)^r
.
\end{equation}
In a similar fashion, we also obtain
%
%
\begin{equation}\label{2.36}
\overline{\xi_{1,2}}(W_+) \le u^\prime\biggl(\frac{\widehat
{L}_0}{\sqrt{d}} \biggr)^{(d-2)} \biggl(4^n \biggl(c
\frac{\widehat{L}_0}{L_0} \biggr)^{(d-2)} \ell_0^{-(n+1)(d-2)} \biggr)^r
.
\end{equation}
The next objective is to show that the trace on $\widehat{C}_2$ of
paths in
the support of $\sum_{1 \le\ell\le r} \rho^\ell_{2,1}$ and $\sum
_{1 \le\ell\le r} \rho^\ell_{1,2}$ is stochastically dominated by
the corresponding trace on $\widehat{C}_2$ of paths in the support of
$\mu
^*_{2,2}$ when $u-u^\prime$ is not too small. An important step is the
following lemma.
\begin{lemma}\label{lem2.2}
For $\ell_0 \ge c(\sqrt{d} + R)$, all $n \ge0$, $m \in I_{n+1}$,
$\cT\in\Lambda_m$, $x\in\partial U$ and $y \in\partial_{\mathrm{int}}
\widehat{C}_2$, one has
%
%
\begin{eqnarray} \label{2.37}
&& P_x[H_{\widehat{C}_1} < R_1 < \infty, X_{R_1} = y]
\nonumber\\[-8pt]\\[-8pt]
&&\qquad\le\biggl(
\frac{\widehat{L}_0}{L_0} \biggr)^{(d-2)} \biggl(
\frac{c}{\ell_0} \biggr)^{(d-2)(n+1)} P_x[H_{\widehat{C}_1} > R_1,
X_{R_1} = y] ,\nonumber
\\
\label{2.38}
&& P_x[H_{\widehat{C}_1} < \infty, R_1 = \infty] \nonumber\\[-8pt]\\[-8pt]
&&\qquad\le\biggl(
\frac{\widehat{L}_0}{L_0} \biggr)^{(d-2)} \biggl(
\frac{c}{\ell_0} \biggr)^{(d-2)(n+1)} P_x[R_1 = \infty= H_{\widehat
{C}_1}] .\nonumber
\end{eqnarray}
\end{lemma}
\begin{pf}
The proof of (\ref{2.37}) closely follows the proof of (\ref{2.30})
in Lemma 2.3 of \cite{SidoSzni09b}. The difference lies in the control
of Harnack constants. Indeed, we first observe that the function $h\dvtx z
\rightarrow P_z[R_1 < \infty$, $X_{R_1} = y] = P_z[H_{\widehat{C}_2} <
\infty$,
$X_{H_{\widehat{C}_2}} = y]$ is a nonnegative function, harmonic in
$\widehat
{C}_2^c$. By (\ref{2.29}), it is therefore harmonic on any $B(z_0,
\frac{L_{n+1}}{20})$ with $z_0 \in\partial U$. One can then find $c$
such that for any $\widetilde{z}, \widetilde{z} ^\prime$ in $\partial
U$, there
exists a sequence $z_i$, $0 \le i \le m$, in $\partial U$ with $m \le
c$, $z_0 = \widetilde{z}$, $z_m = \widetilde{z} ^\prime$ and $|z_{i+1}
- z_i|
\le\frac{L_{n+1}}{100 c_3}$, in the notation of Proposition \ref
{prop1.3}. Indeed, one simply projects $\widetilde{z}, \widetilde{z}
^\prime$
onto the Euclidean sphere in $\IR^d$ of radius $\frac{L_{n+1}}{10}$
with center $\overline{i_2} L_n$, the ``center'' of $U$ [see (\ref
{2.26})] and uses the great circle joining these two points to
construct the sequence.

Using (\ref{1.15}) and a standard chaining argument, it follows that
%
%
\begin{equation}\label{2.39}
\sup_{z \in\partial U} P_z[R_1 < \infty, X_{R_1} = y] \le
c^d \inf_{z \in\partial U} P_z[R_1 < \infty, X_{R_1} = y] .
\end{equation}
The proof of (\ref{2.37}) then proceeds as in Lemma 2.3 of
\cite{SidoSzni09b} [and we use a similar bound to (\ref{2.33}) above, where
$\widehat{C}_1$ replaces $\widehat{C}_2$].

As for (\ref{2.38}), we first note that for $x \in\partial U$, due to
(\ref{1.6}) and (\ref{2.29}), we have
%
%
\begin{eqnarray}\label{2.40}
&&\inf_{x \in\partial U} P_x [R_1 = \infty, H_{\widehat{C}_1}
= \infty] \nonumber\\
&&\qquad\ge1 - 2 2^n \biggl(c \frac{\widehat
{L}_0}{L_{n+1}} \biggr)^{(d-2)} \stackrel{\mbox{\fontsize{8.36}{10.36}\selectfont{(\ref{2.2})}}}{\ge} 1 -
\biggl(\frac{c}{\ell_0} \frac{\widehat
{L}_0}{L_0} \biggr)^{(d-2)} \\
&&\qquad\hspace*{-4.06pt}\stackrel{\mbox{\fontsize{8.36}{10.36}\selectfont{(\ref{2.1})}}}{\ge}
\frac{1}{2},\nonumber
\end{eqnarray}
when $\ell_0 \ge c^\prime(\sqrt{d} + R)$.

On the other hand, a similar calculation leads to
%
%
\begin{eqnarray}\label{2.41}
P_x[H_{\widehat{C}_1} < \infty, R_1 = \infty] & \le & 2^n \biggl(c
\frac{\widehat{L}_{0}}{L_{n+1}} \biggr)^{(d-2)} \nonumber\\[-8pt]\\[-8pt]
&\le&\biggl(
\frac{\widehat{L}_0}{L_0} \biggr)^{(d-2)} \biggl(
\frac{c}{\ell_0} \biggr)^{(d-2)(n+1)}\nonumber
\end{eqnarray}
and (\ref{2.38}) follows.
\end{pf}

The proof of Proposition \ref{prop2.1} then proceeds like the proof of
Proposition 2.3 of \cite{SidoSzni09b} and yields that, under (\ref
{2.19}) (with $u^\prime$ in place of $u_n$ and $u$ in place of $u_{n+1}$),
%
%
\begin{eqnarray}\label{2.42}\hspace*{12pt}
&&\IP[A_{\overline{\cT}_2} (\mu_{2,2})] \nonumber\\
&&\qquad\hspace*{18.21pt} \le
p_n(u^\prime) + \overline{\xi}_{2,1} (W_+) + \overline{\xi}_{1,2} (W_+)
\\
&&\qquad\stackrel{\mbox{\fontsize{8.36}{10.36}\selectfont{(\ref{2.35}), (\ref{2.36})}}}{\le}
p_n(u^\prime) + 2u^\prime\biggl(\frac{\widehat
{L}_0}{\sqrt{d}} \biggr)^{(d-2)} \biggl(4^n \biggl( c
\frac{\widehat{L}_0}{L_0} \biggr)^{(d-2)} \ell_0^{-(n+1)(d-2)} \biggr)^r .\nonumber
\end{eqnarray}
Inserting this inequality into (\ref{2.22}), we thus infer (\ref
{2.20}) under the assumption on (\ref{2.19}).
\end{pf}

We assume from now on that $\ell_0 \ge c(\sqrt{d} + R)$, with $c > 2
c_5$ sufficiently large so that Proposition \ref{prop2.1} holds. We
then choose the sequences $u_n$, $n \ge0$ and $r_n$, $n \ge0$, as follows:
%
%
\begin{eqnarray} \label{2.43}
u_n & = & u_0 \exp\biggl\{ \biggl(\frac{\widehat
{L}_0}{L_0} \biggr)^{(d-2)} \sum_{0 \le k < n} (r_k + 1) \biggl(
\frac{c_5}{\ell_0} \biggr)^{(k +1)(d-2)} \biggr\},
\\
\label{2.44}
r_n & = & r_0 2^n ,
\end{eqnarray}
where $u_0 > 0$ and $r_0$ is a positive integer. The choice (\ref
{2.43}) ensures that (\ref{2.19}) is fulfilled and the increasing
sequence $u_n$ has the finite limit
%
%
\begin{equation}\label{2.45}
u_\infty= u_0 \exp\biggl\{ \biggl(\frac{c_5 \widehat{L}_0}{\ell
_0 L_0} \biggr)^{(d-2)} \biggl(\frac{r_0}{1 - 2 (c_5 \ell
_0^{-1})^{(d-2)}} + \frac{1}{1 - (c_5 \ell
_0^{-1})^{(d-2)}} \biggr) \biggr\}.\hspace*{-35pt}
\end{equation}
The next proposition reduces the task of bounding $p_n(u_n)$ to a set
of conditions which enable us to initiate the induction procedure
suggested by Proposition \ref{prop2.1}. We view $u_\infty$ as a
function of $u_0$, $r_0$, $\ell_0$, $R$ [we introduced $R$ in (\ref{2.1})].
\begin{proposition}\label{prop2.3}
There exists a positive constant $c$ such that when $u_0 > 0$, $r_0 \ge
1$, $\ell_0 \ge c(\sqrt{d} + R)$, $L_0 \ge d$, $\widehat{L}_0 = (\sqrt
{d} + R) L_0$, $R \ge1$ and $K_0 > \log2$ satisfy
%
%
\begin{eqnarray} \label{2.46}
u_\infty\biggl(\frac{\widehat{L}_0}{\sqrt{d}}
\biggr)^{d-2} \vee e^{K_0} &\le& \biggl(\frac{\ell_0
L_0}{c_6 \widehat{L}_0} \biggr)^{{r_0}/{2} (d-2)},
\\ \label{2.47}
p_0 (u_0) &\le& e^{-K_0},
\end{eqnarray}
then
%
%
\begin{equation}\label{2.48}
p_n(u_n) \le e^{-(K_0 - \log2)2^n} \qquad\mbox{for each $n \ge0$} .
\end{equation}
\end{proposition}
\begin{pf}
The argument is similar to Proposition 2.5 of \cite{SidoSzni09b}. We
assume, as mentioned before, that $c > 2c_5$ is large enough so that
Proposition \ref{prop2.1} applies. Condition (\ref{2.46}) implies
that $c_6 \widehat{L}_0 \le\ell_0 L_0$(\mbox{$ =$}$L_1)$. Thus, the last term
in the right-hand side of (\ref{2.20}) satisfies
%
%
\begin{eqnarray}\label{2.49}
&&
u_n \biggl(\frac{\widehat{L}_0}{\sqrt{d}} \biggr)^{(d-2)}
\biggl(4^n \biggl(c_6 \frac{\widehat{L}_0}{L_0}
\biggr)^{(d-2)} \ell_0^{-(n+1)(d-2)} \biggr)^{r_n} \nonumber\\
&&\hspace*{17.6pt}\qquad\le
u_\infty\biggl(\frac{\widehat{L}_0}{\sqrt{d}}
\biggr)^{(d-2)}\biggl (c_6 \frac{\widehat{L}_0}{ \ell_0
L_0} \biggr)^{(d-2)r_n} \biggl( \frac{4}{\ell
_0^{d-2}} \biggr)^{nr_n} \\
&&\qquad\stackrel{\mbox{\fontsize{8.36}{10.36}
\selectfont{(\ref{2.2}), (\ref{2.46})}}}{\le}
\biggl(c_6 \frac{\widehat{L}_0}{\ell_0 L_0}
\biggr)^{{r_n}/{2} (d-2)} .\nonumber
\end{eqnarray}
As a result, (\ref{2.20}) yields that for $n \ge0$,
%
%
\begin{equation}\label{2.50}
p_{n+1}(u_{n+1}) \le p_n(u_n) \biggl(p_n(u_n) + \biggl(c_6
\frac{\widehat{L}_0}{\ell_0 L_0} \biggr)^{({r_0}/{2})
2^n(d-2)} \biggr).
\end{equation}
We then define by induction $K_n$, $n \ge0$, via the following
relation valid for $n \ge1$:
%
%
\begin{equation}\label{2.51}\qquad\quad
K_n = K_0 - \sum_{0 \le n^\prime< n} 2^{-(n^\prime+ 1)} \log\biggl(
1 + e^{K_{n^\prime} 2^{n^\prime}} \biggl(c_6 \frac{
\widehat{L}_0}{\ell_0 L_0} \biggr)^{({r_0}/{2}) 2^{n^\prime}
(d-2)} \biggr)
\end{equation}
so that $K_n \le K_0$ and hence
%
%
\begin{eqnarray}\label{2.52}\quad
K_n & \ge & K_0 - \sum_{n^\prime\ge0} 2^{-(n^\prime+ 1)} \log\biggl(
1 + e^{K_0 2^{n^\prime}} \biggl(c_6 \frac{ \widehat
{L}_0}{\ell_0 L_0} \biggr)^{({r_0}/{2}) 2^{n^\prime}(d-2)} \biggr)
\nonumber\\[-8pt]\\[-8pt]
& \stackrel{\mbox{\fontsize{8.36}{10.36}\selectfont{(\ref{2.46})}}}{\ge} & K_0 - \sum_{n^\prime\ge0}
2^{-(n^\prime+ 1)} \log2 = K_0 - \log2 >
0 .\nonumber
\end{eqnarray}
As we now show by induction, we have $p_n(u_n) \le e^{-K_n 2^n}$.

Indeed, this inequality holds for $n = 0$, due to (\ref{2.47}), and if
it holds for $n \ge0$, then, due to (\ref{2.50}), we find that
\begin{eqnarray*}
p_{n+1}(u_{n+1}) & \le & e^{-K_n 2^n} \biggl(e ^{-K_n 2^n} + \biggl(c_6
\frac{\widehat{L}_0}{\ell_0 L_0} \biggr)^{
({r_0}/{2}) 2^n(d-2)} \biggr)
\\
& = & e^{-K_n 2^{n+1}} \biggl(1 + e^{K_n 2^n} \biggl(c_6
\frac{\widehat{L}_0}{\ell_0 L_0} \biggr)^{({r_0}/{2})
2^n(d-2)} \biggr) \\
&\stackrel{\mbox{\fontsize{8.36}{10.36}\selectfont{(\ref{2.51})}}}{=}&
e^{-K_{n+1} 2^{n+1}}.
\end{eqnarray*}
This proves that $p_n(u_n) \le e^{-K_n 2^n}$ for all $n \ge0$ and
(\ref{2.48}) follows.
\end{pf}
\begin{remark}\label{rem2.4}
One of the main issues we now have to face is proving the local
estimate $p_0(u_0) \le e^{-K_0}$ [see (\ref{2.47})] for large $d$,
with $u_0$ of order close to $\log d$ (and a posteriori close to
$u_*$). We further need $K_0$ sufficiently large so that $2^{-(K_0 -
\log2)2^n}$ overcomes the combinatorial complexity arising from the choice
of the binary trees in the upper bound (\ref{2.18}), that is, overcomes
the\vspace*{2pt} factor $|\Lambda_n|
\stackrel{\mbox{\fontsize{8.36}{10.36}\selectfont{(\ref{2.13})}}}{\le} (c_4 \ell
_0)^{2(d-1)(2^n-1)}$. Devising this local estimate will be the aim of
the next section and will involve aspects of random interlacements at a
shorter range, where features reminiscent of random interlacements on
$2d$-regular trees (cf. Section 5 of \cite{Teix09b}) will be evident.
\end{remark}

\section{Local connectivity bounds}\label{sec3}

The aim of this section is to derive exponential bounds on the decay of
the probability of existence of a path in the vacant set at level $u_0
= (1 +5 \varepsilon) \log d$, starting at the origin and traveling at
$\ell^1$-distance $Md$, where $M$ is an arbitrary integer and $d \ge
c(\varepsilon, M)$ (cf. Corollary \ref{cor3.4}). For this purpose, we
develop an enhanced Peierls-type argument. The main step appears in
Theorem \ref{theo3.1} below. In the present section, aspects of random
interlacements on $\IZ^d$ for large $d$, reminiscent of random
interlacements on $2d$-regular trees (cf. \cite{Teix09b}) will play a
important role. We introduce the parameter
%
%
\begin{equation}\label{3.1}
0 < \varepsilon< \tfrac{1}{3}.
\end{equation}
We also introduce, in the notation of (\ref{1.14}),
%
%
\begin{equation}\label{3.2}
L= c_7 d\qquad \mbox{where } c_7 = [e^8 c_2] + 2 .
\end{equation}
The main result of this section is the following estimate on the
connectivity function.
\begin{theorem}[($d \ge c$)]\label{theo3.1}
For any positive integer $M$, we have
%
%
\begin{equation}\label{3.3}\quad
\IP\bigl[0 \stackrel{\cV^{u_0}}{\longleftrightarrow} S_1 (0,ML)\bigr] \le
\exp\biggl\{ \frac{M(M-1)}{2} L + 3 Md -
\frac{\varepsilon^2}{5} Md \log d \biggr\},
\end{equation}
where the notation is similar to (\ref{2.14}) and
%
%
\begin{equation}\label{3.4}
u_0 = ( 1 + 5 \varepsilon) \log d .
\end{equation}
\end{theorem}
\begin{pf}
Observe that any self-avoiding path from $0$ to $S_1 (0,ML)$
successively visits the $\ell^1$-spheres $S_1(0,iL)$, $i = 0, \ldots,
M-1$. Thus, considering the first $[\frac{\varepsilon}{10} d]$ steps
of the path consecutive to the successive entrances in the various
spheres $S_1(0,iL)$, we obtain $M$ self-avoiding paths $\pi_i$, $i =
0, \ldots, M-1$, where $\pi_i$ starts in $S_1(0,iL)$ and has $[\frac
{\varepsilon}{10} d]$ steps for each $i$. Denoting by $z_i$, $i =
0,\ldots, M-1$, the respective starting points of these paths, we find that
%
%
\begin{equation}\label{3.5}\quad
\IP\bigl[0 \stackrel{\cV^{u_0}}{\longleftrightarrow} S_1 (0,ML)\bigr] \le
\sum_{z_i,\pi_i} \IP[\cV^{u_0} \supseteq\operatorname{range} \pi_i
\mbox{ for } i = 0, \ldots, M-1] ,
\end{equation}
where the above sum runs over $z_i \in S(0,iL)$ and self-avoiding paths
$\pi_i$ with $[\frac{\varepsilon}{10} d]$ steps and starting points
$z_i$, $i = 0,\ldots, M-1$. The next lemma provides a very rough bound
on the cardinality of $\ell^1$-spheres and $\ell^1$-balls. Crucially,
it shows that $\ell^1$-spheres and balls of radius $cd$ are ``rather
small,'' that is, their cardinality grows at most geometrically in $d$.
\begin{lemma}[($\ell\in\IN$)]\label{lem3.2}
%
%
\begin{eqnarray}\label{3.6}
&&\mbox{\hphantom{i}\textup{(i)}} \quad |S_1 (0,\ell) | \le2^d e^{\ell+ d},
\nonumber\\[-8pt]\\[-8pt]
&&\mbox{\textup{(ii)}} \quad |B_1(0,\ell) | \le2^d e^{\ell+ 1 + d} .\nonumber
\end{eqnarray}
\end{lemma}
\begin{pf}
We express the generating function of $|S_1(0,k)|$, $k \ge0$, as
follows. Given $|t| < 1$, we have

%
\begin{eqnarray}\label{3.7}
\sum_{k \ge0} t^k |S_1(0,k)| & = & \sum_{k \ge0} t^k \mathop{\sum
_{m_1,\ldots,m_d \ge0}}_{m_1 + \cdots+ m_d = k} 2^{|\{i \in\{
1,\ldots,d\}; m_i \not= 0\}|}
\nonumber\\
& = &\sum_{m_1,\ldots,m_d \ge0} t^{m_1 + \cdots+ m_d} 2^{|\{i \in\{
1,\ldots, d\}; m_i \not= 0\}|} \nonumber\\[-8pt]\\[-8pt]
& = &\biggl(1 + 2 \sum_{m \ge1} t^m \biggr)^d
\nonumber\\
& = &\biggl( \frac{1+ t}{1-t} \biggr)^d \le\frac{2^d}{(1-t)^d}.\nonumber
\end{eqnarray}
As a result, we see that for $0 < t < 1$, $\ell\ge0$,
\[
|S_1(0,\ell)| \le2^d(1-t)^{-d} t^{-\ell} .
\]
Choosing $t = \ell/ (d+\ell)$, we find that
%
%
\begin{equation}\label{3.8}
|S_1 (0,\ell)| \le2^d \biggl(1 + \frac{\ell
}{d} \biggr)^d \biggl(1 + \frac{d}{\ell} \biggr)^\ell
\le2^d e^{\ell+ d} ,
\end{equation}
where we have used the inequality $1 + u \le e^u$ in the last step.
This proves (\ref{3.6})(i). As for the inequality (\ref{3.6})(ii), by
(\ref{3.6})(i), we can write
%
%
\begin{equation}\label{3.9}
|B_1(0,\ell)| \le2^d e^d \sum^\ell_{k=0} e^k = 2^d e^d
\frac{e^{\ell+1}-1}{e-1} \le2^d e^{\ell+ 1 + d}
\end{equation}
and our claims follows.
\end{pf}

We now come back to (\ref{3.5}). By a very rough counting argument for
the number of possible choices of $\pi_i$, we have a Peierls-type bound:
%
%
\begin{eqnarray}\label{3.10}
&&\IP\bigl[0 \stackrel{\cV^{u_0}}{\longleftrightarrow} S_1(0,ML)\bigr]\nonumber \\
&&\hspace*{7.75pt}\qquad\le
\Biggl(\prod^{M-1}_{k=0} |S_1(0,k L)| (2 d)^{
{\varepsilon}/{10} d} \Biggr) \nonumber\\
&&\hspace*{7.75pt}\qquad\quad{}\times\sup_{z_i,\pi_i} \IP[\cV
^{u_0} \supseteq\operatorname{range} \pi_i, i = 0, \ldots, M-1]
\nonumber\\
&&\qquad\stackrel{\mbox{\fontsize{8.36}{10.36}\selectfont{(\ref{3.6})(i)}}}{\le}
\Biggl(\prod^{M-1}_{k=0} 2^d e^{k L+d} \Biggr) (2d)^{
{\varepsilon}/{10} M d} \\
&&\hspace*{7.75pt}\qquad\quad{}\times\sup_{z_i,\pi_i} \IP[\cV^{u_0}
\supseteq\operatorname{range} \pi_i, i = 0, \ldots, M-1] \nonumber\\
&&\hspace*{7.75pt}\qquad\le
e^{{M(M-1)}/{2} L + 2 M d} (2d)^{{\varepsilon}/{10} M d}\nonumber\\
&&\hspace*{7.75pt}\qquad\quad{}\times
\sup_{z_i,\pi_i} \IP[\cV^{u_0} \supseteq\operatorname{range}
\pi_i, i = 0, \ldots, M-1] ,\nonumber
\end{eqnarray}
where the supremum runs over a similar collection as the sum in (\ref{3.5}).

The next objective is to bound the probability in the last line of
(\ref{3.10}). For this purpose, for each $x$ in the set
%
%
\begin{eqnarray}\label{3.11}\qquad
\hspace*{-12pt}B \stackrel{\mathrm{def}}{=} \bigcup^{M-1}_{i=0} B_1 \biggl(z_i,
\frac{\varepsilon}{10} d \biggr) \nonumber\\[-8pt]\\[-8pt]
\eqntext{\mbox
{(pairwise disjoint $\ell^1$-balls appear in this union)},}
\end{eqnarray}
we write $z_x$ for the unique $z_i$ such that $x \in B (z_i,
\frac{\varepsilon}{10} d )$. We then define, for any
$x$ in $B$, the subset $W^*_x$ of $W^*$---see above (\ref
{1.24})---(not to be confused with $W^*_{\{x\}}$):
%
%
\begin{eqnarray} \label{3.12}
W_x^* & = & \mbox{the image under $\pi^*$ of}\nonumber\\
&&{}  \biggl\{\mbox{$w \in W$: the
minimum of $d_1(z_x,w(n))$, $n \in\IZ$,}
\nonumber\\[-8pt]\\[-8pt]
&&\hspace*{7.4pt} \mbox{is reached for the first time at $w(n) = x$ and $w$ }\nonumber\\
&&\hspace*{11pt} \mbox{does not
enter any }B_1 \biggl(z_i, \frac{\varepsilon}{10} d
\biggr) \mbox{ with } z_i \not= z_x \biggr\} .\nonumber
\end{eqnarray}
Note that, clearly, $W^*_x \subseteq W^*_{\{x\}}$ and that
%
%
\begin{equation}\label{3.13}
\mbox{$W^*_x$, $x \in B$, are pairwise disjoint measurable subsets of
$W^*$} .
\end{equation}
It then follows that for $z_i, \pi_i$, $0 \le i \le M-1$, as in (\ref
{3.10}), we have
%
%
\begin{eqnarray}\label{3.14}
&&\IP[\cV^{u_0} \supseteq\operatorname{range} \pi_i, i = 0,\ldots, M-1]
\nonumber\\
&&\qquad\le
\IP\Biggl[\omega\Biggl(\bigcup^{M-1}_{i=0} \bigcup
_{x \in\operatorname{range} \pi_i} W_x^* \times[0,u_0] \Biggr) = 0
\Biggr]\nonumber\\[-8pt]\\[-8pt]
&&\qquad =
\exp\Biggl\{- u_0 \sum^{M-1}_{i=0} \sum_{x \in\operatorname{range} \pi
_i} \nu(W^*_x) \Biggr\}
\nonumber\\
&&\qquad\le\exp\biggl\{ - u_0 M \frac{\varepsilon d}{10}
\times\inf_{x \in B} \nu(W^*_x) \biggr\} .\nonumber
\end{eqnarray}
We will now seek a lower bound on $\nu(W^*_x)$ for $x \in B$.

Choosing $K = \{x\}$ in (\ref{1.25}), (\ref{1.26}), by (\ref{1.19}),
we see that for any $x$ in $B$,
%
%
\begin{eqnarray}\label{3.15}
\nu(W^*_x) & = & P_x \bigl[|X_n - z_x|_1 \ge|x-z_x|_1\mbox{, for $n
\ge0$,}\nonumber\\
&&\hspace*{39.9pt}\mbox{and } H_{\bigcup_{z_i \not= z_x} B_1(z_i,
{\varepsilon}/{10} d)} = \infty\bigr]
\nonumber\\
&&{} \times P_x \bigl[|X_n - z_x|_1> |x-z_x|_1\mbox{, for $n >
0$,}\nonumber\\[-8pt]\\[-8pt]
&&\hspace*{53.2pt} \mbox{and }
H_{\bigcup_{z_i \not= z_x} B_1(z_i, {\varepsilon
}/{10} d)} = \infty\bigr]
\nonumber\\
& \ge & \biggl(P_x [|X_n - z_x|_1 > |x-z_x|_1\mbox{,  for } n > 0]\nonumber\\
&&\hspace*{38.2pt}{} -
\sum_{z_i \not= z_x} P_x
\bigl[ H_{B_1(z_i, {\varepsilon}/{10}d)} < \infty\bigr] \biggr)^2_+.\nonumber
\end{eqnarray}
In view of (\ref{1.14}) and the choice of $L$ in (\ref{3.2}), we see
that when $d \ge8$, we have
%
%
\begin{eqnarray}\label{3.16}
&&\sum_{z_i \not= z_x} P_x \bigl[ H_{B_1(z_i, {\varepsilon
}/{10} d)} < \infty\bigr] \nonumber\\
&&\qquad\hspace*{21.28pt}\le\sum_{z_i \not= z_x} \sup
_{y \in B(z_i, {\varepsilon}/{10} d)} g(y-x) \biggl|B_1
\biggl(0, \frac{\varepsilon}{10} d \biggr) \biggr|\nonumber\\[-8pt]\\[-8pt]
&&\qquad\stackrel{\mbox{\fontsize{8.36}{10.36}\selectfont{(\ref{1.14}), (\ref{3.6})(ii)}}}{\le}
2 \sum_{j \ge1} \biggl(\frac{c_2 d}{j L -
{\varepsilon}/{5} d} \biggr)^{{d}/{2} - 2} 2^d e^{
{\varepsilon}/{10} d + 1 + d} \nonumber\\
&&\hspace*{17.1pt}\qquad\stackrel{\mbox{\fontsize{8.36}{10.36}\selectfont{(\ref{3.2})}}}{\le}
2e^{-8({d}/{2} - 2) + 3d + 1} \sum_{j \ge1} j^{-({d}/{2}
- 2)} \stackrel{d \ge8}{\le} c e^{-d} .\nonumber
\end{eqnarray}
The next lemma yields a lower bound on the first term in the last line
of (\ref{3.15}).
\begin{lemma}[$(d \ge c)$]
When $|y|_1 \le\frac{d}{2}$, one has
%
%
\begin{equation}\label{3.17}
P_y [ |X_n |_1 > |y|_1 \mbox{ for all } n > 0] \ge1 -
\frac{4(|y|_1 \vee1)}{2d - (|y|_1 \vee1)} .
\end{equation}
\end{lemma}
\begin{pf}
We first note that for $z = (z_1,\ldots,z_d)$ in $\IZ^d$, $P_z$-a.s.,
$ | |X_1|_1 - |z|_1 | = 1$, and
%
%
\begin{eqnarray}\label{3.18}
P_z[|X_1|_1 = |z|_1 + 1] &=& \frac{1}{2d} \Biggl(2d -
\sum^d_{k=1} 1\{z_k \not= 0\} \Biggr) \ge
p_{|z|_1}\nonumber\\[-8pt]\\[-8pt]
\eqntext{\mbox{where }
\displaystyle p_m \stackrel{\mathrm{def}}{=} \biggl(\frac{1}{2} +
\frac{1}{2} \biggl( 1 - \frac{m}{d} \biggr)_+\biggr)
\qquad\mbox{for $m \ge0$} .}
\end{eqnarray}
We then introduce the canonical Markov chain $N_n$ on $\IN$ that jumps
to $m+1$ with probability $p_m$ and to $m-1$ with probability $q_m =
1-p_m$ when located at $m$. We denote by $Q_m$ the canonical law of
this Markov chain starting in $m$. In view of (\ref{3.18}), a coupling
argument shows that we can construct $X_n$ and $N_n$ on the same
probability space so that a.s. $|X_n|_1 \ge N_n$ for all\vspace*{1pt} $n \ge0$ and
$X_0 = y \in\IZ^d$, $N_0 = |y|_1$. Consequently, we see that when $y
\not= 0$, we have the bound (with $m = |y|_1 \le\frac{d}{2}$)
%
%
\begin{eqnarray}\label{3.19}
&&
P_y \bigl[ H_{S_1(0,d^2)} < \widetilde{H}_{B_1(0,|y|_1)} \bigr] \nonumber\\
&&\qquad\ge Q_{|y|_1} \bigl[H_{d^2}
< \widetilde{H}_{|y|_1}\bigr]
\\
&&\qquad =  p_m (1 + \rho_{m+1} + \rho_{m+1} \rho_{m+2} + \cdots+ \rho
_{m+1} \cdots\rho_{d^2 -
1})^{-1},\nonumber
\end{eqnarray}
where $\rho_\ell= \frac{q_\ell}{p_\ell}$ for $\ell\ge0$ and we
have used \cite{Chun60}, (5), page 73.

Note that the expression in the right-hand side of (\ref{3.19}) is a
decreasing function of each $\rho_\ell$, $m + 1 < \ell< d^2$. If we
further observe that $\rho_\ell\le(\frac{1}{2} - \frac{1}{2}
\times\frac{1}{4})(\frac{1}{2} + \frac{1}{2} \times\frac
{1}{4})^{-1} = \frac{3}{5}$ for $m + 1 < \ell\le\frac{3}{4} d$
and $\rho_\ell\le1$ for $\frac{3}{4} d < \ell\le d^2 - 1$, then
we see that the above expression is bigger than
\begin{eqnarray*}
&&\biggl(1 - \frac{m}{2d} \biggr) \biggl( 1 +
\frac{m}{2d-m} \sum_{k \ge0} \biggl(
\frac{3}{5} \biggr)^k + \biggl(
\frac{3}{5} \biggr)^{[{3}/{4} d]-m} d^2 \biggr)^{-1}\\
&&\hspace*{5.41pt}\qquad\stackrel{m \le{d/2}}{\ge}
\biggl(1 - \frac{m}{2d} \biggr) \biggl( 1 +
\frac{5}{2} \frac{m}{2d-m} +
\frac{5}{3} \biggl(
\frac{3}{5} \biggr)^{{d}/{4}} d^2 \biggr)^{-1} \\
&&\qquad\mathop{\ge
}_{1 \le
m \le{d/2}}^{d \ge c} \biggl(1 -
\frac{m}{2d} \biggr) \biggl( 1 + 3
\frac{m}{2d-m} \biggr)^{-1}
\\
&&\hspace*{13.94pt}\qquad\ge\biggl(1 - \frac{m}{2d} \biggr) \biggl( 1 - 3
\frac{m}{2d-m} \biggr) \qquad\biggl( \mbox{$\ge$}0 \mbox{ since } m \le
\frac{d}{2} \biggr) .
\end{eqnarray*}
By the strong Markov property at time $H_{S_1(0,d^2)}$, we thus find
that for $d \ge c$, $1 \le|y|_1 \le\frac{d}{2}$, we have
%
%
\begin{eqnarray}\label{3.20}\qquad
&&P_y [|X_n|_1 > |y|_1 \mbox{ for all } n > 0] \nonumber\\
&&\hspace*{4.3pt}\qquad\ge
\biggl(1 - \frac{|y|_1}{2d} \biggr) \biggl(1 -
\frac{3|y|_1}{2d - |y|_1} \biggr) \Bigl(1 - \sup
_{|z|_1= d^2} P_z\bigl[H_{B_1(0,{d}/{2})}< \infty\bigr] \Bigr)\\
&&\qquad \stackrel
{\mbox{\fontsize{8.36}{10.36}\selectfont{(\ref{1.1})}}}{\ge}
\biggl(1 - \frac{|y|_1}{2d} \biggr) \biggl(1 -
\frac{3|y|_1}{2d - |y|_1} \biggr) - \sup_{|z|\ge
d^{3/2}} P_z \bigl[H_{B(0,d)} < \infty\bigr] \nonumber\\
&&\qquad\stackrel{\mbox{\fontsize{8.36}{10.36}\selectfont{(\ref{1.6})}}}{\ge}
1 - \frac{|y|_1}{2d} - \frac
{3|y|_1}{2d - |y|_1} + \frac{3|y|^2_1}{2d(2d -
|y|_1)} - \biggl(\frac{c}{\sqrt{d}} \biggr)^{(d-2)}\nonumber\\
&&\qquad\mathop{\ge}_{y \not= 0}^{d \ge c} 1 - \frac
{4|y|_1}{2d - |y|_1} .\nonumber
\end{eqnarray}
This completes the proof of (\ref{3.17}) when $y \not= 0$. The
extension to the case $y=0$ is immediate.
\end{pf}
\noqed\end{pf}

We use the above lemma to bound the first term in the last line of
(\ref{3.15}) from below. In view of (\ref{3.16}) and (\ref{3.17}),
we thus find that for $d \ge c$ and any $x \in B$ [see (\ref{3.11})],
%
%
\begin{equation}\label{3.21}\qquad
\nu(W^*_x) \ge\biggl(1 - 5 \frac{|x-z_x|_1 \vee
1}{2d-(|x-z_x|_1 \vee1)} \biggr)^2 \ge1 - 10 \frac
{\varepsilon/10}{2-\varepsilon/10} \ge1 - \varepsilon.
\end{equation}
Coming back to (\ref{3.14}), we thus find that
%
%
\begin{equation}\label{3.22}
\IP[\cV^{u_0} \supseteq\operatorname{range} \pi_i, i = 0,\ldots,M-1] \le
\exp\biggl\{ - \frac{u_0}{10} M \varepsilon(1 -
\varepsilon) d \biggr\}.
\end{equation}
Inserting this bound into the last line of (\ref{3.10}), we obtain
\begin{eqnarray*}
&&\IP\bigl[0 \stackrel{\cV^{u_0}}{\longleftrightarrow} S_1(0,ML)\bigr] \\
&&\hspace*{4pt}\qquad\le
\exp\biggl\{ \frac{M(M-1)}{2} L+2 Md -
\frac{u_0}{10} \varepsilon(1-\varepsilon) Md \biggr\}
(2d)^{{\varepsilon}/{10} Md}
\\
&&\qquad\stackrel{\mbox{\fontsize{8.36}{10.36}\selectfont{(\ref{3.4})}}}{\le}
\exp\biggl\{ \frac{M(M-1)}{2} L+3 Md + \frac{\varepsilon}{10}
Md \log d\\
&&\qquad\quad\hspace*{57.2pt} - \frac{1}{10} (\varepsilon+ 4
\varepsilon^2 - 5 \varepsilon^3) Md \log d \biggr\}.
\end{eqnarray*}
Since $5\varepsilon^3 \le2\varepsilon^2$, due to (\ref{3.1}), the
claim (\ref{3.3}) follows.

We will use the following corollary in the proof of Theorem \ref
{theo0.1} in the next section.
\begin{corollary}[{[with (\ref{3.1}), (\ref{3.2})]}]\label{cor3.4}
If $M \ge1$, then for $d \ge c(M,\varepsilon)$,
%
%
\begin{equation}\label{3.23}
\IP\bigl[0 \stackrel{\cV^{u_0}}{\longleftrightarrow} S_1(0,ML)\bigr] \le
\exp\biggl\{ - \frac{\varepsilon^2}{10} dM \log
d \biggr\}.
\end{equation}
\end{corollary}
\begin{pf}
This is an immediate consequence of (\ref{3.3}).
\end{pf}
\begin{remark}\label{rem3.5}
One should note that the bound of Theorem \ref{theo3.1} deteriorates
when $M$ becomes large. One can view Theorem \ref{theo3.1} as a
Peierls-type bound (slightly enhanced due to the role of $M$ in the
proof). In the next section, we will choose $M$ as a large constant
depending on $\varepsilon$ and use Corollary \ref{cor3.4} to produce
the local estimate which will enable us to initiate the renormalization
scheme of Section \ref{sec2}. In this way, the local estimate on crossings in
$\cV^{u_0}$ at $\ell^1$-distance of order $c(\varepsilon)d$ will be
transformed into an estimate on crossings at all scales in $\cV
^{u_\infty}$, where $u_\infty\le(1 + 10 \varepsilon) \log d$.
\end{remark}

\section{Denouement}\label{sec4}

In this section, we prove Theorem \ref{theo0.1}. We combine the local
bound on the connectivity function at level $u_0$ of the last section
(cf. Corollary \ref{cor3.4}) with the renormalization scheme of
Section \ref{sec2} (cf. Proposition \ref{prop2.3}) in order to produce a bound
on vacant crossings at a level $u_\infty\in[(1 + 5 \varepsilon) \log
d$, $(1
+ 10 \varepsilon) \log d]$, valid at arbitrarily large scales.
\begin{pf*}{Proof of Theorem \ref{theo0.1}}
We choose $\varepsilon$ and $u_0$
as in
(\ref{3.1}), (\ref{3.4}), respectively. For the renormalization
scheme of Section \ref{sec2}, we choose [the constant $c_7$ appears in (\ref{3.2})]
%
%
\begin{equation} \label{4.1}
L_0 = d,\qquad \widehat{L}_0 = \bigl(\sqrt{d} + R\bigr) L_0\qquad \mbox{with $R = 300
c_7 \varepsilon^{-2}$}
\end{equation}
and
\begin{equation}
\label{4.2}
\ell_0 = d .
\end{equation}
In the notation of Proposition \ref{prop2.3} and (\ref{2.13}), we choose
%
%
\begin{equation} \label{4.3}
r_0 = 24
\end{equation}
and
\begin{equation}
\label{4.4}
K_0 = \log\bigl(4 (c_4 \ell_0)^{2(d-1)}\bigr)
\stackrel{\mbox{\fontsize{8.36}{10.36}\selectfont{(\ref{4.2})}}}{=}
\log\bigl(4 (c_4 d)^{2(d-1)}\bigr) .
\end{equation}
In the application of Corollary \ref{cor3.4}, we choose
%
%
\begin{equation}\label{4.5}
M = [100 \varepsilon^{-2}] + 1
\end{equation}
so that in the notation of (\ref{3.2}), (\ref{4.1}),
%
%
\begin{equation}\label{4.6}
ML + 1 \le R L_0 .
\end{equation}
We will now check that the assumptions of Proposition \ref{prop2.3}
hold for $d \ge c(\varepsilon)$. By (\ref{2.45}), we see that for $d \ge
c(\varepsilon)$,
%
%
\begin{equation}\label{4.7}
u_0 = (1 + 5 \varepsilon) \log d < u_\infty< ( 1+ 10 \varepsilon) \log d
\end{equation}
and also that
%
%
\begin{equation}\label{4.8}
\widehat{L}_0 \le2d^{{3/2}} .
\end{equation}
As a result, we find that
%
%
\begin{equation}\label{4.9}
u_\infty\biggl(\frac{\widehat{L}_0}{\sqrt{d}}
\biggr)^{(d-2)} \le(1 + 10 \varepsilon) (\log d) (2d)^{(d-2)}
\end{equation}
and that
%
%
\begin{equation}\label{4.10}
e^{K_0} = 4(c_4 d)^{2(d-1)} ,
\end{equation}
whereas, on the other hand,
%
%
\begin{equation}\label{4.11}
\biggl(\frac{\ell_0 L_0}{c_6 \widehat{L}_0}
\biggr)^{{r_0}/{2} (d-2)} \ge(cd)^{6(d-2)} .
\end{equation}
Since $2(d-1) < 6(d-2)$, we see that for $d \ge c(\varepsilon)$, the
expression in the left-hand side of (\ref{4.11}) dominates the
corresponding expressions in (\ref{4.9}) and (\ref{4.10}), that is,
(\ref{2.46}) holds.

There remains to check (\ref{2.47}). For this purpose, we apply
Corollary \ref{cor3.4} and find that for $d \ge c(\varepsilon)$, since
$\widehat
{L}_0 \ge\sqrt{d} L_0 + M L + 1$ [cf. (\ref{4.1}), (\ref{4.6})],
we have
%
%
\begin{eqnarray}\label{4.12}
p_0(u_0) & = & \IP\bigl[ [0,L_0 - 1]^d \stackrel{\cV
^{u_0}}{\longleftrightarrow} \partial_{\mathrm{int}} B(0,\widehat{L}_0) \bigr]
\nonumber\\
& \le & L^d_0 \IP\bigl[0 \stackrel{\cV^{u_0}}{\longleftrightarrow} S_1(0,ML)\bigr]
\mathop{\le}_{\mbox{\fontsize{8.36}{10.36}\selectfont{(\ref{4.5})}}}
^{\mbox{\fontsize{8.36}{10.36}\selectfont{(\ref{3.23})}}} \exp\{ d \log d - 10 d
\log d\} \\
&=&
d^{-9d}.\nonumber
\end{eqnarray}
We thus find that for $d \ge c(\varepsilon)$, $p_0(u_0) \le e^{-K_0}$, that
is, (\ref{2.47}) holds as well. It now follows from Proposition \ref
{prop2.3} that for $d \ge c(\varepsilon)$,
%
%
\begin{equation}\label{4.13}
p_n (u_\infty) \le e^{-(K_0 - \log2)2^n} \qquad\mbox{for all $n \ge0$} .
\end{equation}
Taking (\ref{2.13}), (\ref{2.18}) into account yields that for all $n
\ge1$,
%
%
\begin{eqnarray}\label{4.14}
&&\IP\bigl[[0,L_n - 1]^d \stackrel{\cV^{u_\infty
}}{\longleftrightarrow} \partial_{\mathrm{int}} [-L_n, 2L_n - 1]^d \bigr]
\nonumber\\[-8pt]\\[-8pt]
&&\qquad\le(c_4 \ell_0)^{2(d-1)(2^n-1)} e^{-(K_0 - \log2)2^n}
\stackrel{\mbox{\fontsize{8.36}{10.36}\selectfont{(\ref
{4.10})}}}{\le} 2^{-2^n} .\nonumber
\end{eqnarray}
In particular, the above inequality implies that $\IP[0 \stackrel{\cV
^{u_\infty}}{\longleftrightarrow} \infty] = 0$ and hence $u_* \le
u_\infty< (1 + 10 \varepsilon) \log d$ for $d \ge c(\varepsilon)$. The
claim (\ref
{0.6}) readily follows. Combining this upper bound with the lower bound
(\ref{0.3}), we have thus proven Theorem \ref{theo0.1}.
\end{pf*}
\begin{remark}\label{rem4.1}

(1) The inequality (\ref{4.14}), together with the fact that $L_n =
L_0 \ell_0^n$ for $n \ge0$, is more than enough to show that for
$\varepsilon$ as in (\ref{3.1}) and $d \ge c(\varepsilon)$,
\[
\lim_{L \rightarrow\infty} L^\gamma\IP\bigl[B_\infty(0,L) \stackrel
{\cV^{(1+10 \varepsilon)\log d}}{\longleftrightarrow} S_\infty(0,2L)\bigr] =
0 ,
\]
for some, and, in fact, all, $\gamma> 0$. From the definition of the
critical parameter $u_{**}$ in \cite{Szni09c},
%
%
\begin{eqnarray}\label{4.15}
u_{**} & = &\inf\{u \ge0; \alpha(u) > 0\}\hspace*{150pt}
\nonumber\\[-8pt]\\[-8pt]
\eqntext{\mbox{where }\displaystyle\alpha(u) = \sup\Bigl\{\alpha\ge0; \lim_{L \rightarrow\infty}
L^\alpha\IP\bigl[B_\infty(0,L) \stackrel{\cV^u}{\longleftrightarrow}
S_\infty(0,2L)\bigr] =
0\Bigr\} }
\end{eqnarray}
(the supremum is, by convention, equal to zero when the set is empty),
we thus find that for $d \ge c(\varepsilon)$,
%
%
\begin{equation}\label{4.16}
u_{**} \le( 1 + 10 \varepsilon) \log d .
\end{equation}
Since $u_* \le u_{**}$, it follows that we have also proven that
%
%
\begin{equation}\label{4.17}
\lim_d u_{**} \big/ \log d = 1 .
\end{equation}
It is presently an open question whether $u_* = u_{**}$; however, we
know that $0 < u_* \le u_{**} < \infty$ for all $d \ge3$ (cf.
\cite{Szni09e}) and that for $u > u_{**}$, the connectivity function has a
stretched exponential decay (cf. \cite{SidoSzni09b}).

(2) One may wonder whether the following reinforcement of (\ref{0.4})
actually holds:
\[
\IP[0 \in\cV^{u_*}] = e^{-u_*/g(0)} \sim(2d)^{-1} \qquad\mbox{as } d
\rightarrow\infty.
\]
This would indicate a similar high-dimensional behavior as for
Bernoulli percolation; see
\cite{AlonBenjStac04,BollKoha94,Gord91,HaraSlad90,Kest90}. In the case of
interlacement percolation on a $2d$-regular tree, such an asymptotic
behavior is known to hold (cf. \cite{Teix09a}).
\end{remark}

\begin{appendix}\label{app}
\section*{Appendix}

In this appendix, we prove an elementary inequality which is involved
in the proof of the Green function estimate (\ref{1.14}); see Lemma
\ref{A.1} below. We then prove, in Lemma \ref{lemA.2}, a bound on
Harnack constants in terms of killed Green functions for
nearest-neighbor Markov chains on graphs. The result is stated in a
rather general formulation due to the fact that it is of independent
interest. It is an adaptation of Lemma 10.2 of \cite{GrigTelc01}. We
recall that Lemma \ref{lemA.2} enters the proof of Proposition \ref{prop1.3}.
\setcounter{lem}{0}
\begin{lem}\label{lemA.1}
%
%
\setcounter{equation}{0}
\begin{eqnarray}\label{A.1}
&&\mbox{For } a,b \ge0\qquad \sqrt{a^2 + b^2} \log\bigl(1 + \sqrt{a^2
+ b^2}\bigr) \nonumber\\[-8pt]\\[-8pt]
&&\hspace*{52.8pt}\qquad\qquad\le a \log(1+a) + b \log(1+b) .\nonumber
\end{eqnarray}
\end{lem}
\begin{pf}
We introduce $\psi(u) = u \log(1 + u)$, $u \ge0$, as well as
$\varphi_b(a) = \sqrt{a^2 + b^2}$ and $\chi_b(a) = \psi(a) + \psi
(b) - \psi(\varphi_b(a))$ for $a,b \ge0$. We want to show that
%
%
\begin{equation}\label{A.2}
\chi_b(a) \ge0 \qquad\mbox{for } a,b > 0 .
\end{equation}
We note that $\chi_b(0) = 0$ and that
\[
\chi^\prime_b(a) = \log(1 + a) + 1 - \frac
{1}{1+a} - \biggl(\log\bigl(1 + \varphi_b(a)\bigr) + 1 - \frac
{1}{1 + \varphi_b(a)} \biggr) \frac{a}{\varphi
_b(a)} .
\]
The claim (\ref{A.2}) will follow once we show that
%
%
\begin{equation}\label{A.3}
\chi^\prime_b(a) \ge0\qquad \mbox{for $a,b > 0$} .
\end{equation}
To this end, we note that for $a > 0$, $\chi^\prime_0(a) = 0$ and that
%
%
\begin{eqnarray}\label{A.4}
\frac{\partial}{\partial b} \chi^\prime_b(a) &
= &- \biggl( \frac{1}{1+\varphi_b(a)}
\frac{b}{\varphi_b(a)} + \frac{1}{(1+\varphi
_b(a))^2} \frac{b}{\varphi_b(a)} \biggr)
\frac{a}{\varphi_b(a)}
\nonumber\\
&&{} + \biggl(\log\bigl(1 + \varphi_b(a)\bigr) + 1 - \frac
{1}{1+\varphi_b(a)} \biggr) \frac{a b}{\varphi_b(a)^3}
\nonumber\\[-8pt]\\[-8pt]
& = &\frac{a b}{(1 + \varphi_b(a))\varphi_b(a)^3}\nonumber\\
&&\times{}
\biggl\{ \log\bigl(1 + \varphi_b(a)\bigr) \bigl(1 + \varphi_b(a)\bigr) -
\frac{\varphi_b(a)}{1+\varphi_b(a)} \biggr\} .\nonumber
\end{eqnarray}
We introduce the function $\rho(u) = \log(1 + u) (1 + u) - \frac
{u}{1+u}$, $u \ge0$. Observe that $\rho(0) = 0$ and $\rho^\prime(u)
= \log(1 + u) + 1 - \frac{1}{(1+u)^2} \ge0$ so that $\rho(u) \ge0$
for $u \ge0$. Coming back to the last line of (\ref{A.4}), we find
that for $a > 0$, $\frac{\partial}{\partial b} \chi^\prime_b(a)
\ge0$ for $b \ge0$. This shows (\ref{A.3}) and the claim (\ref
{A.1}) then follows.
\end{pf}

We now turn to the second result of this appendix. We consider a
connected graph $\Gamma$ with an at most countable vertex set $E$ and
edge set $\mathcal{E}$ (a subset of the collection of unordered pairs
of $E$).
Given $U \subseteq E$, we define $\partial U$, $\partial_{\mathrm{int}} U$
and $\overline{U}$ similarly to what is described at the beginning of the
Section \ref{sec1} (with obvious modifications). We consider an irreducible
Markov chain on $E$, nearest-neighbor in the broad sense (i.e., at each
step, the Markov chain moves to a vertex which is at graph-distance at
most $1$ from its current location). We write $X_n$, $n \ge0$, for the
canonical process, $P_x$ for the canonical law starting from $x \in E$
and otherwise use similar notation as described at the beginning of
Section \ref{sec1}. We denote by $p(x,y)$, $x,y \in E$, the transition
probability. We assume that the Markov chain satisfies the following
ellipticity condition:
%
%
\begin{equation}\label{A.5}
p(x,y) > 0 \qquad\mbox{when $x,y$ are neighbors (i.e., $\{x,y\} \in\mathcal{E}$)}.
\end{equation}
For $f$ a bounded function on $E$, we define
%
%
\begin{equation}\label{A.6}
L f(x) = E_x[f(X_1)] - f(x) = \sum_{y \sim x} p(x,y)\bigl(f(y) - f(x)\bigr)
\qquad\mbox{for $x \in E$} ,\hspace*{-35pt}
\end{equation}
where $y \sim x$ means that $y = x$ or $y$ is a neighbor of $x$. Given
$U \subseteq E$, a bounded function on $\overline{U}$ is said to be
\textit{harmonic in} $U$ when (with a slight abuse of notation)
%
%
\begin{equation}\label{A.7}
L f(x) = 0 \qquad\mbox{for $x \in U$} .
\end{equation}
When $U$ is a finite strict subset of $E$, the Green function killed
outside $U$ is defined as follows (the notation is similar to that in
Section \ref{sec1}):
%
%
\begin{equation}\label{A.8}
G_U(x,y) = E_x \biggl[\sum_{k \ge0} 1\{X_k = y, T_U > k\} \biggr],\qquad x,y
\in E .
\end{equation}
It follows from the ellipticity assumption (\ref{A.5}) that when $U$
is connected, $G_U(x,y) > 0$ for all $x,y \in U$. The next lemma is an
adaptation of Lemma 10.2 of~\cite{GrigTelc01}.
\begin{lemma}\label{lemA.2}
Assume that $\varnothing\not= U_1 \subseteq U_2 \subseteq U_3$ are
finite strict subsets of $E$, with $U_3$ connected, and that $u$ is a
bounded nonnegative function on $\overline{U}_3$ which is harmonic in
$U_3$. We then have
%
%
\begin{equation}\label{A.9}
\max_{U_1} u \le K \min_{U_1} u ,
\end{equation}
where
%
%
\begin{equation}\label{A.10}
K = \max_{x,y \in U_1} \max_{z \in\partial_{\mathrm{int}} U_2} G_{U_3} (x,z) / G_{U_3}(y,z) .
\end{equation}
\end{lemma}
\begin{pf}
We define, for $x \in E$,
%
%
\begin{equation}\label{A.11}
v(x) = E_x [u (X_{H_{U_2}}), H_{U_2} < T_{U_3}] .
\end{equation}
We first note that
%
%
\begin{equation}\label{A.12}\quad
u(x) \ge v(x) \qquad\mbox{for $x \in\overline{U}_3$}\quad\mbox{and}\quad u(x) = v(x)\qquad\mbox{for $x
\in U_2$} .
\end{equation}
Indeed, in view of (\ref{A.11}), $u$ and $v$ agree on $U_2$ and,
thanks to our assumptions, $u(X_{n \wedge T_{U_3}})$, $n \ge0$, is a
bounded martingale under $P_x$, $x \in\overline{U}_3$, so that by the
stopping theorem, we find that
\begin{eqnarray*}
u(x) &=& E_x [u(X_{H_{U_2} \wedge T_{U_3}})]  = v(x) + E_x
[u(X_{T_{U_3}}), T_{U_3} < H_{U_2}]
\\
&\ge& v(x) \qquad\mbox{for } x \in\overline{U}_3 .
\end{eqnarray*}
The claim (\ref{A.12}) then follows.

Applying the simple Markov property at time 1 in (\ref{A.11}), when $x
\in U_3 \setminus U_2$, we see that
%
%
\begin{equation}\label{A.13}
\mbox{$v$ is harmonic in $U_3 \setminus U_2$} .
\end{equation}
In addition, we have, for $x \in U_2$,
\[
v(x) = u(x) = \sum_{y \sim x} p(x,y) u(y) \stackrel{\mbox{\fontsize{8.36}{10.36}\selectfont{(\ref
{A.12})}}}{\ge} \sum_{y \sim x} p(x,y) v(y)
\]
and the last inequality is an equality when $x \in U_2 \setminus
\partial_{\mathrm{int}} U_2$. We have thus shown that
%
%
\begin{equation}\label{A.14}
Lv = 1_{\partial_{\mathrm{int}} U_2} Lv \le0 \qquad\mbox{on $U_3$} .
\end{equation}
Applying the stopping theorem, we see that, under any $P_x$,
\[
v(X_{n \wedge T_{U_3}}) - \sum_{0 \le k < n \wedge T_{U_3}}
Lv(X_k),\qquad
n \ge0,\qquad \mbox{is a martingale} .
\]
Taking expectations and letting $n$ tend to infinity, we obtain the identity
%
%
\begin{eqnarray}\label{A.15}
v(x) & = & E_x[v(X_{T_{U_3}})] - E_x \biggl[\sum_{0 \le k < T_{U_3}} Lv
(X_k) \biggr]
\nonumber\\
& = & - \sum_{z \in E} G_{U_3} (x,z) Lv(z) \\
&\stackrel{\mbox{\fontsize{8.36}{10.36}\selectfont{(\ref{A.14})}}}{=}&
\sum_{z \in\partial_{\mathrm{int}} U_2} G_{U_3} (x,z)
(-Lv)(z),\qquad x \in
E .\nonumber
\end{eqnarray}
Since $v$ and $u$ agree on $U_2 \supseteq U_1$, (\ref{A.9}) is a
direct consequence of the above representation formula for $v$.
\end{pf}
\end{appendix}

%

%
\printaddresses


\begin{thebibliography}{23}

\bibitem{AlonBenjStac04}
%
\begin{barticle}[mr]
\bauthor{\bsnm{Alon},~\bfnm{Noga}\binits{N.}},
\bauthor{\bsnm{Benjamini},~\bfnm{Itai}\binits{I.}} \AND
\bauthor{\bsnm{Stacey},~\bfnm{Alan}\binits{A.}}
(\byear{2004}).
\btitle{Percolation on finite graphs and isoperimetric inequalities}.
\bjournal{Ann. Probab.}
\bvolume{32}
\bpages{1727--1745}.
\bid{doi={10.1214/009117904000000414}, mr={2073175}}
\end{barticle}
%
\endbibitem

\bibitem{BollKoha94}
%
\begin{barticle}[mr]
\bauthor{\bsnm{Bollob{\'a}s},~\bfnm{B.}\binits{B.}} \AND
\bauthor{\bsnm{Kohayakawa},~\bfnm{Y.}\binits{Y.}}
(\byear{1994}).
\btitle{Percolation in high dimensions}.
\bjournal{European J. Combin.}
\bvolume{15}
\bpages{113--125}.
\bid{doi={10.1006/eujc.1994.1014}, mr={1261058}}
\end{barticle}
%
\endbibitem

\bibitem{BroaHamm57}
%
\begin{barticle}[mr]
\bauthor{\bsnm{Broadbent},~\bfnm{S.~R.}\binits{S.~R.}} \AND
\bauthor{\bsnm{Hammersley},~\bfnm{J.~M.}\binits{J.~M.}}
(\byear{1957}).
\btitle{Percolation processes. {I}. {C}rystals and mazes}.
\bjournal{Proc. Cambridge Philos. Soc.}
\bvolume{53}
\bpages{629--641}.
\bid{mr={0091567}}
\end{barticle}
%
\endbibitem

\bibitem{CernTeixWind09}
%
\begin{bmisc}[vtex]
\bauthor{\bsnm{\v{C}ern\'y},~\bfnm{J.}\binits{J.}},
\bauthor{\bsnm{Teixeira},~\bfnm{A.}\binits{A.}} \AND
\bauthor{\bsnm{Windisch},~\bfnm{D.}\binits{D.}}
(\byear{2009}).
\bhowpublished{Giant vacant component left by a random walk in a random
$d$-regular graph. Preprint. Available at}
\url{http://www.math.ethz.ch/\textasciitilde cerny/publications.html}.
\end{bmisc}
%
\endbibitem

\bibitem{Chun60}
%
\begin{bbook}[vtex]
\bauthor{\bsnm{Chung},~\bfnm{Kai~Lai}\binits{K.~L.}}
(\byear{1960}).
\btitle{Markov Chains with Stationary Transition Probabilities}.
\bpublisher{Springer}, \baddress{Berlin}.
\bid{mr={0116388}}
\end{bbook}
%
\endbibitem

\bibitem{Gord91}
%
\begin{barticle}[mr]
\bauthor{\bsnm{Gordon},~\bfnm{Daniel~M.}\binits{D.~M.}}
(\byear{1991}).
\btitle{Percolation in high dimensions}.
\bjournal{J. London Math. Soc. (2)}
\bvolume{44}
\bpages{373--384}.
\bid{doi={10.1112/jlms/s2-44.2.373}, mr={1136447}}
\end{barticle}
%
\endbibitem

\bibitem{GrigTelc01}
%
\begin{barticle}[mr]
\bauthor{\bsnm{Grigor'yan},~\bfnm{Alexander}\binits{A.}} \AND
\bauthor{\bsnm{Telcs},~\bfnm{Andras}\binits{A.}}
(\byear{2001}).
\btitle{Sub-{G}aussian estimates of heat kernels on infinite graphs}.
\bjournal{Duke Math. J.}
\bvolume{109}
\bpages{451--510}.
\bid{doi={10.1215/S0012-7094-01-10932-0}, mr={1853353}}
\end{barticle}
%
\endbibitem

\bibitem{Grim99}
%
\begin{bbook}[mr]
\bauthor{\bsnm{Grimmett},~\bfnm{Geoffrey}\binits{G.}}
(\byear{1999}).
\btitle{Percolation},
\bedition{2nd} ed.
\bseries{Grundlehren der Mathematischen Wissenschaften [Fundamental Principles
of Mathematical Sciences]}
\bvolume{321}.
\bpublisher{Springer}, \baddress{Berlin}.
\bid{mr={1707339}}
\end{bbook}
%
\endbibitem

\bibitem{HaraSlad90}
%
\begin{barticle}[mr]
\bauthor{\bsnm{Hara},~\bfnm{Takashi}\binits{T.}} \AND
\bauthor{\bsnm{Slade},~\bfnm{Gordon}\binits{G.}}
(\byear{1990}).
\btitle{Mean-field critical behaviour for percolation in high dimensions}.
\bjournal{Comm. Math. Phys.}
\bvolume{128}
\bpages{333--391}.
\bid{mr={1043524}}
\end{barticle}
%
\endbibitem

\bibitem{HaraSlad92}
%
\begin{barticle}[mr]
\bauthor{\bsnm{Hara},~\bfnm{Takashi}\binits{T.}} \AND
\bauthor{\bsnm{Slade},~\bfnm{Gordon}\binits{G.}}
(\byear{1992}).
\btitle{The lace expansion for self-avoiding walk in five or more dimensions}.
\bjournal{Rev. Math. Phys.}
\bvolume{4}
\bpages{235--327}.
\bid{doi={10.1142/S0129055X9200008X}, mr={1174248}}
\end{barticle}
%
\endbibitem

\bibitem{Kest90}
%
\begin{bincollection}[vtex]
\bauthor{\bsnm{Kesten},~\bfnm{Harry}\binits{H.}}
(\byear{1990}).
\btitle{Asymptotics in high dimensions for percolation}.
In \bbooktitle{Disorder in Physical Systems}
\bpages{219--240}.
\bpublisher{Oxford Univ. Press}, \baddress{New York}.
\bid{mr={1064563}}
\end{bincollection}
%
\endbibitem

\bibitem{Lawl91}
%
\begin{bbook}[vtex]
\bauthor{\bsnm{Lawler},~\bfnm{Gregory~F.}\binits{G.~F.}}
(\byear{1991}).
\btitle{Intersections of Random Walks}.
\bpublisher{Birkh\"auser}, \baddress{Boston, MA}.
\bid{mr={1117680}}%
\end{bbook}%
%
\endbibitem%

\bibitem{Mont56}
%
\begin{barticle}[mr]
\bauthor{\bsnm{Montroll},~\bfnm{Elliot~W.}\binits{E.~W.}}
(\byear{1956}).
\btitle{Random walks in multidimensional spaces, especially on periodic
lattices}.
\bjournal{J. Soc. Indust. Appl. Math.}
\bvolume{4}
\bpages{241--260}.
\bid{mr={0088110}}
\end{barticle}
%
\endbibitem

\bibitem{Olve74}
%
\begin{bbook}[vtex]
\bauthor{\bsnm{Olver},~\bfnm{F.~W.~J.}\binits{F.~W.~J.}}
(\byear{1974}).
\btitle{Asymptotics and Special Functions}.
\bpublisher{Academic Press}, \baddress{New York}.
\bid{mr={0435697}}
\end{bbook}
%
\endbibitem

\bibitem{Pang93}
%
\begin{barticle}[mr]
\bauthor{\bsnm{Pang},~\bfnm{M.~M.~H.}\binits{M.~M.~H.}}
(\byear{1993}).
\btitle{Heat kernels of graphs}.
\bjournal{J. London Math. Soc. (2)}
\bvolume{47}
\bpages{50--64}.
\bid{doi={10.1112/jlms/s2-47.1.50}, mr={1200977}}
\end{barticle}
%
\endbibitem

\bibitem{SidoSzni09a}
%
\begin{barticle}[mr]
\bauthor{\bsnm{Sidoravicius},~\bfnm{Vladas}\binits{V.}} \AND
\bauthor{\bsnm{Sznitman},~\bfnm{Alain-Sol}\binits{A.-S.}}
(\byear{2009}).
\btitle{Percolation for the vacant set of random interlacements}.
\bjournal{Comm. Pure Appl. Math.}
\bvolume{62}
\bpages{831--858}.
\bid{doi={10.1002/cpa.20267}, mr={2512613}}
\end{barticle}
%
\endbibitem

\bibitem{SidoSzni09b}
%
\begin{barticle}[vtex]
\bauthor{\bsnm{Sidoravicius},~\bfnm{V.}\binits{V.}} \AND
\bauthor{\bsnm{Sznitman},~\bfnm{A.~S.}\binits{A.~S.}}
(\byear{2010}).
\btitle{Connectivity bounds for the vacant set of random interlacements}.
\bjournal{Ann. Inst. H. Poincar\'e Probab. Statist.}
\bvolume{46}
\bpages{976--990}.
\end{barticle}
%
\endbibitem

\bibitem{Szni07a}
%
\begin{barticle}[auto:SpringerTagBib|2009-01-14|16:51:27]
\bauthor{\bsnm{Sznitman},~\bfnm{A.~S.}\binits{A.~S.}}
(\byear{2010}).
\btitle{Vacant set of random interlacements and percolation}.
\bjournal{Ann. Math.}
\bvolume{171}
\bpages{2039--2087}.
\end{barticle}
%
\endbibitem

\bibitem{Szni09c}
%
\begin{barticle}[mr]
\bauthor{\bsnm{Sznitman},~\bfnm{Alain-Sol}\binits{A.-S.}}
(\byear{2009}).
\btitle{Upper bound on the disconnection time of discrete cylinders and random
interlacements}.
\bjournal{Ann. Probab.}
\bvolume{37}
\bpages{1715--1746}.
\bid{doi={10.1214/09-AOP450}, mr={2561432}}
\end{barticle}
%
\endbibitem

\bibitem{Szni09d}
%
\begin{barticle}[vtex]
\bauthor{\bsnm{Sznitman},~\bfnm{Alain-Sol}\binits{A.-S.}}
(\byear{2009}).
\btitle{On the domination of random walk on a discrete cylinder by random
interlacements}.
\bjournal{Electron. J. Probab.}
\bvolume{14}
\bpages{1670--1704}.
\bid{mr={2525107}}
\end{barticle}
%
\endbibitem

\bibitem{Szni09e}
%
\begin{bmisc}[vtex]
\bauthor{\bsnm{Sznitman},~\bfnm{A.~S.}\binits{A.~S.}}
(\byear{2010}).
\bhowpublished{A lower bound on the critical parameter of interlacement
percolation in high dimension.
\textit{Probab. Theory Related Fields}. To appear}.
Available at \href{http://arxiv.org/abs/1003.0334}{arXiv:1003.0334}.
\end{bmisc}
%
\endbibitem

\bibitem{Teix09a}
%
\begin{barticle}[mr]
\bauthor{\bsnm{Teixeira},~\bfnm{Augusto}\binits{A.}}
(\byear{2009}).
\btitle{On the uniqueness of the infinite cluster of the vacant set of random
interlacements}.
\bjournal{Ann. Appl. Probab.}
\bvolume{19}
\bpages{454--466}.
\bid{doi={10.1214/08-AAP547}, mr={2498684}}
\end{barticle}
%
\endbibitem

\bibitem{Teix09b}
%
\begin{barticle}[vtex]
\bauthor{\bsnm{Teixeira},~\bfnm{A.}\binits{A.}}
(\byear{2009}).
\btitle{Interlacement percolation on transient weighted graphs}.
\bjournal{Electron. J. Probab.}
\bvolume{14}
\bpages{1604--1628}.
\bid{mr={2525105}}
\end{barticle}
%
\endbibitem

\end{thebibliography}
\end{document}